%% file: main.tex
\documentclass{amsart}
\input{preamble}
\input{macros}

\calclayout

\title{Two Categorifications of the Local Langlands Correspondence for Tori}
\author{Ruide Fu}

\begin{document}
\begin{abstract}
The stack of local Langlands parameters for a torus is a Picard stack.
In this article, we explicitly determine its Picard dual and show that
the Fourier--Mukai transform gives rise to the integral categorical
local Langlands correspondence for the torus. This is the
categorification of the local Langlands correspondence and answers a
conjecture of X.~Zhu. Moreover, we establish a geometric version of
this correspondence, whose categorical trace reproduces the previous
result.
\end{abstract}

\maketitle
\tableofcontents

\input{content-v4}

{\footnotesize
\printbibliography
}
\end{document}

%% file: preamble.tex
\usepackage{amsmath}
\usepackage{amsthm}
\usepackage{amssymb}
\usepackage[backend=biber,style=alphabetic]{biblatex}
\usepackage{mathrsfs}
\usepackage{mathtools}
\usepackage{verbatim}
\usepackage{dcounter}
\usepackage{geometry}
\usepackage{tikz-cd}
\usepackage{tikz}
\usepackage[inline]{enumitem}
\usepackage{xcolor}
\usepackage{graphicx}
\usepackage{stackrel}
\usepackage{hyperref}
\hypersetup{colorlinks=false,
  linkcolor=red,
  citecolor=green,
  urlcolor=blue,
  pdfborder={0 0 1}
}

\numberwithin{equation}{section}

\newcounter{point}[section]
\renewcommand{\thepoint}{\thesection.\arabic{point}}
\newcommand{\makepoint}[1][]{
  \ifnum \value{point}>0 \medskip \fi
  \stepcounter{point}
  \noindent
  \thepoint. \textbf{#1}}
\newcommand{\refmakepoint}[1][]{
  \ifnum \value{point}>0 \medskip \fi
  \refstepcounter{point}
  \noindent
  \thepoint. \textbf{#1}
}

\theoremstyle{definition}
\newtheorem{definition}[point]{Definition}
\newtheorem{rmk}[point]{Remark}
\newtheorem{remark}[point]{Remark}

\theoremstyle{plain}
\newtheorem{thm}[point]{Theorem}
\newtheorem{conj}[point]{Conjecture}
\newtheorem{lem}[point]{Lemma}
\newtheorem{prop}[point]{Proposition}
\newtheorem{cor}[point]{Corollary}
\newtheorem{eg}[point]{Example}

\newtheorem{lemma}[point]{Lemma}

\makeatletter
\renewcommand{\thesubsection}{
  \ifnum \c@subsection<1 \@arabic\c@section
  \else \thesection.\@arabic\c@subsection
  \fi
}
\makeatother

\newcommand{\nhd}{neighbourhood {}}

\newcommand{\wrt}{with respect to {}}
\newcommand{\parti}{particular {}}

\DeclareMathOperator{\rank}{rank}
\DeclareMathOperator{\Gal}{Gal}
\DeclareMathOperator{\Tr}{Tr}

\renewcommand{\Im}{\operatorname{Im}}

\DeclareMathOperator{\Spec}{Spec}
\DeclareMathOperator{\Ad}{Ad}

\DeclareMathOperator{\id}{id}
\DeclareMathOperator{\Hom}{Hom}
\DeclareMathOperator{\Ext}{Ext}
\DeclareMathOperator{\Tor}{Tor}

\DeclareMathOperator{\Res}{Res}
\DeclareMathOperator{\res}{res}
\DeclareMathOperator{\Ind}{Ind}
\DeclareMathOperator{\ind}{ind}
\DeclareMathOperator{\cores}{cor}
\DeclareMathOperator{\cofib}{cofib}

\DeclareMathOperator{\colim}{colim}

\DeclareMathOperator{\Coh}{Coh}

\DeclareMathOperator{\QCoh}{QCoh}
\DeclareMathOperator{\IndCoh}{IndCoh}
\DeclareMathOperator{\Bun}{Bun}

\newcommand{\Gm}{\mathbb{G}_m}
\newcommand{\BGm}{\mathbb{BG}_m}

\newcommand{\m}{\mathfrak{m}}

\newcommand{\Z}{\mathbb{Z}}

\newcommand{\F}{\mathbb{F}}
\newcommand{\B}{\mathbb{B}}

\newcommand{\inj}{\hookrightarrow}

\newcommand{\vsim}{\rotatebox[origin=c]{270}{$\sim$}}

\newcommand{\bb}[1]{\mathbb{#1}}      

\newcommand{\mcal}[1]{\mathcal{#1}}   

\newcommand{\mbf}[1]{\mathbf{#1}}     


%% file: macros.tex
\usepackage{pifont}

\DeclareMathOperator{\Homm}{\underline{\Hom}}
\newcommand{\cHom}{\Hom_{\mathrm{cts}}}

\newcommand{\scHom}{\Hom_{\mathrm{sc}}}
\DeclareMathOperator{\Shv}{Shv}
\newcommand{\Sing}{\operatorname{Sing}}
\newcommand{\Nilp}{0}

\newcommand{\rev}{\mathrm{rev}}
\DeclareMathOperator{\HH}{HH}
\newcommand{\Lincat}{\mathrm{Lincat}}
\newcommand{\Lincatlam}{\Lincat_\Lambda}

\newcommand{\Loc}{\mathrm{Loc}}
\newcommand{\Locc}{\Loc_{\cT, F}}
\newcommand{\Locn}{\Locc^{(n)}}
\newcommand{\uLoc}{\underline{\Loc}}
\newcommand{\n}{{\geq n}}
\newcommand{\np}{{(n)}}
\newcommand{\uLocn}{\uLoc^{(n)}}
\newcommand{\Y}{{\Loc^\square_{\cT, \breve F}}}
\newcommand{\X}{{\Loc_{\cT, \breve F}}}
\DeclareMathOperator{\Lsig}{\mathcal L_\sigma}
\DeclareMathOperator{\Spf}{Spf}

\newcommand{\W}{W_{\F/F}}
\newcommand{\Wn}{W^{(n)}}
\renewcommand{\F}{E}

\newcommand{\Ex}{E^\times}
\newcommand{\Rx}{R^\times}

\renewcommand{\Tor}{\boldsymbol{\mathrm{Tor}}_{T, \mathrm{iso}_F}}
\newcommand{\TorS}{\boldsymbol{\mathrm{Tor}}_{S, \mathrm{iso}_F}}
\newcommand{\TorTE}{\boldsymbol{\mathrm{Tor}}_{T, \mathrm{iso}_E}}

\newcommand{\BT}{\mathfrak{B}(T)}
\newcommand{\isoc}{\mathrm{Isoc}}

\newcommand{\T}{\mathcal{TOR}}
\newcommand{\LnT}{L^\n T}
\newcommand{\LNT}{\LnT}
\newcommand{\doublequotientLT}{\raisebox{-0.7ex}{$L^\n T$} \backslash \raisebox{0.7ex}{$LT$} / \raisebox{-0.7ex}{$L^\n T$}}

\newcommand{\BLT}{\B LT}

\newcommand{\utf}{\mathscr T}
\newcommand{\uto}{\mathscr T_0}
\newcommand{\BTF}{\B T(F)}
\newcommand{\BTO}{\B T(F)_0}

\newcommand{\Rep}{\operatorname{Rep}}

\newcommand{\Ch}{\operatorname{Ch}}
\newcommand{\ch}{\operatorname{ch}}
\newcommand{\alg}{\mathrm{alg}}
\newcommand{\Zalg}{Z^1_\alg}
\newcommand{\Halg}{H^1_\alg}
\newcommand{\repfun}{\mathcal R}
\newcommand{\repfunc}{\mathcal R^c}
\newcommand{\repfunsc}{\mathcal R^{sc}}
\newcommand{\rfc}[1]{\repfunc_{#1}}
\newcommand{\rfsc}[1]{\repfunsc_{#1}}

\newcommand{\cT}{{}^L T}

\newcommand{\I}{I^{\mathrm{rel}}_{E/F}}
\newcommand{\ab}{\mathrm{ab}}
\newcommand{\Iab}{I_E^\ab}

\newcommand{\piLT}{\pi_1(L^+T)}
\newcommand{\piLTvalue}{(\hat L\otimes\Iab)^I}
\newcommand{\tame}{\mathrm{tame}}

\newcommand{\Fl}{{\mathbb{F}_\ell}}

\newcommand{\Zl}{{\mathbb{Z}_\ell}}
\newcommand{\Flb}{\overline{\mathbb{F}}_\ell}
\newcommand{\Qlb}{\overline{\mathbb{Q}}_\ell}
\newcommand{\Zlb}{\overline{\mathbb{Z}}_\ell}
\newcommand{\Adsig}{\Ad_\sigma}
\newcommand{\lx}{\Lambda^\times}

\newcommand{\tphi}{\tilde{\phi}}
\newcommand{\tpsi}{\tilde{\psi}}
\newcommand{\hH}{\mathscr{H}}
\newcommand{\hZ}{\mathscr{Z}}

\newcommand{\shD}{\mathcal{D}}
\newcommand{\Corr}{\mathrm{Corr}}
\newcommand{\perf}{\mathrm{perf}}
\newcommand{\PreStk}{\mathrm{PreStk}}
\newcommand{\AlgStk}{\mathrm{AlgStk}}
\newcommand{\IndArStk}{\mathrm{IndArStk}}
\newcommand{\cH}{\mathrm{H}}
\newcommand{\cV}{\mathrm{V}}
\newcommand{\HR}{\mathrm{HR}}
\newcommand{\VR}{\mathrm{VR}}
\newcommand{\All}{\mathrm{All}}
\newcommand{\IndEproet}{\mathrm{IndEproet}}
\newcommand{\Sch}{\mathrm{Sch}}

%% file: content-v4.tex

\section{Introduction}
\label{introduction}
\noindent
Let $G$ be a reductive group over a non-archimedean local field $F$.
X.~Zhu carried out a detailed study of the stack of local and global
Langlands parameters $\Loc_{^c G, F}$ defined over $\Z[1/p]$, and
proposed a categorical form of the arithmetic local Langlands
correspondence \cite{Zhu21}. The automorphic side of this Langlands
correspondence is given by the category of $\ell$-adic sheaves on the
stack of isocrystals $\isoc_G = LG/\mathrm{Ad}_\sigma LG$, where
$\sigma$ is the Frobenius on the residue field $\kappa_F$ of $F$. A
central part of the conjecture is as follows.

\begin{conj}[Zhu21, Conjecture 4.6.4]
    \label{llc}
	Let $(G, B, T, e)$ be a pinned quasi-split reductive group over $F$
	and $\Lambda = \Zlb, \Qlb$ or $\Flb$. Fix an additive character
	$\psi_0: F\to \Lambda^\times$ with conductor $\mathcal O_F$. Then
	there is a natural equivalence of stable $\infty$-categories
	\[
		\mathbb{L}_G:\Shv(\isoc_G, \Lambda) \to
		\IndCoh_{\hat{\mathcal N}_{^c G}}(\Loc_{^c G, F})
	\]
	sending the Whittaker sheaf to the structural sheaf
	$\mathcal{O}_{\Loc_{^c G, F}}$. The category
	\[
        \IndCoh_{\hat{\mathcal N}_{^c G}}(\Loc_{^c G, F})
        \subset \IndCoh(\Loc_{^c G, F})
    \]
    is a subcategory of ind-coherent sheaves with certain singular
    support condition.
\end{conj}

In a more recent development, Zhu established the tame categorical local
Langlands correspondence with $\Qlb$-coefficients and the unipotent
correspondence with $\Flb$-coefficients \cite{Zhu25}.

\makepoint \textbf{The First Categorification.}
In the first part of this paper, we verify Conjecture \ref{llc} in the
case $G = T$, where $T$ is an arbitrary tori over $F$. In this case,
no choice of a Whittaker datum $\psi_0$ is required.

A key feature of the torus case is that $\Locc$ can be defined
over $\Z$, and there exists an analogue of $\isoc_T$ also
defined over $\Z$. To explain the latter, recall that
the stack $\isoc_T$ is a disjoint union of copies
of the classifying stack $[*/T(F)]$ indexed by the Kottwitz set
$B(T)$. We denote by $\Tor$ the algebraic stack over $\Z$ that
also consists of copies of $[*/T(F)]$ indexed by $B(T)$.
We say $\Tor$ is an analogue of $\isoc_T$ because for
$\Lambda = \Zlb, \Qlb$, or $\Flb$, we have
\begin{equation}
    \label{shv-qcoh}
	\Shv(\isoc_T, \,\Lambda)\cong\QCoh(\Tor\otimes\Lambda).
\end{equation}

The key ingredient of the proof is another conjecture of Zhu.

\begin{conj}[Zhu21, Conjecture 3.2.2]
        There is a natural Poincar\'e line bundle $\mathcal L_T$ on
        $\Tor\times\Locc$ inducing an isomorphism
        $\Locc^\vee\cong\Tor$.
\end{conj}

We define the Poincar\'e line bundle over $\Z$, and show the
isomorphism $\Locc^\vee\cong\Tor$ when $2$ is invertible
in the base. We shall see that
the isomorphism is essentially taking the cup product with the
fundamental class from local class field theory.
The proof of this conjecture is largely inspired by the work of
R.~Langlands \cite{Lan97}. 

The Fourier--Mukai transform was first introduced to the Langlands
program in the work of Laumon, where the transform is applied to
construct the geometric Langlands correspondence for $GL(1)$. Later,
in the work of Braverman-Bezrukavnikov \cite{BB07} and Chen-Zhu
\cite{CZ14}, Fourier--Mukai is used to establish a generic version of
the Langlands correspondence in positive characteristic (for any
reductive group $G$). In the same spirit, we use the Fourier--Mukai
transform to establish the categorical local Langlands correspondence
for tori with integral coefficients:

\begin{thm}
    \label{llt}
    Let $p$ be invertible in $\Lambda$.
    The Fourier--Mukai transform via the Poincar\'e line bundle
	$\mathcal L_T$ gives the equivalence of stable $\infty$-categories
	\[
		\bb{L}: \IndCoh_{\Nilp}(\Locc\otimes\Lambda)
            \cong\QCoh(\Tor\otimes\Lambda).
	\]
\end{thm}

\makepoint \textbf{The Second Categorification.}
Let $\Lambda = \Flb$ or $\Qlb$. Let $\mathfrak T$ be the
connected N\'eron model of $T$, and let $L^+T$ be the positive loop
group. There is an isomorphism between the abelian group of continuous
characters of Serre's fundamental group $\pi_1(L^+T)$ and the abelian
group of character sheaves on $L^+T$:
\[
	\cHom(\pi_1(L^+T), \lx) = 
	\mathrm{CS}(L^+T, \Lambda)
.\]

The second main result of this paper categorifies this
isomorphism:
\begin{thm}
    There exists a fully-faithful, $t$-exact, monoidal functor
    \[
        \Ch: \IndCoh(\X) \to \Shv(LT, \Lambda),
    \]
    where $\X$ is the stack of inertial Langlands parameters
    of $F$ in the dual group $\cT$. The essential image of $\Ch$ is
    compactly generated by the translations of all character sheaves
    on $L^+T$ to $LT$.
\end{thm}

This theorem can be regarded as a ``second-level categorification'' in
the following sense. Both $\IndCoh(\X)$ and $\Shv(LT, \Lambda)$ carry
canonical Frobenius structures that are intertwined by $\Ch$. The
decategorification of $\Ch$ via the categorical trace construction is
precisely the categorical local Langlands correspondence $\bb{L}$:

\begin{prop}
    There is a commutative diagram
    \begin{equation}
    \begin{tikzcd}
        \Tr(\IndCoh(\X), \sigma)
        \ar[r] \ar[d, "\sim", sloped]
        & \Tr(\Shv(LT), \sigma) \ar[d, "\sim", sloped] \\
        \IndCoh_\Nilp(\Locc) \ar[r, "\bb{L}"] & \Shv(\isoc_T),
    \end{tikzcd}
    \end{equation}
    where $\bb{L}$ is the categorical local Langlands correspondece
    under the identification \(
    \Shv(\isoc_T, \,\Lambda)\cong\QCoh(\Tor\otimes\Lambda)
    \).
    Furthermore, the two vertical arrows are canonical equivalences. 
\end{prop}

Some of the results in this paper were independently obtained by
K.~Zou \cite{Zou24} by different methods. There are two main
differences between this work and Zou's. First, Zou considers integral
$\ell$-adic sheaves on $\Bun_G$ over the Fargues-Fontaine curve, while
we consider quasi-coherent sheaves on $\Tor$, which can be integral.
Second, to establish the Langlands correspondence, Zou uses the
spectral action on the Whittaker sheaf, while we use the Fourier--Mukai
transform.

The outline of the paper is as follows. Section~\ref{notations}
introduces the definitions of $\Locc$ and $\Tor$, along with notations
for the first half of the paper.
Section~\ref{short-exact-sequence} establishes an important split
exact sequence involving $\Locc$.
Section~\ref{duality} constructs the Poincar\'e line bundles and
proves the first main theorem.
Section~\ref{transform} discusses the Fourier--Mukai transform
and the categorical local Langlands correspondence for tori.
Section~\ref{notation-ii} introduces the stack of inertial Langlands
parameters $\X$, the loop group $LT$ and all other notations for the
second half of the paper.
Section~\ref{geom-llc} proves the second main theorem and the
result of decategorification.

\textbf{Acknowledgment.}
I am deeply grateful to my advisor Xinwen Zhu for suggesting this
problem and for his generous guidance and support throughout this
work. I also thank Justin Campbell, Tamir Hemo, and Longke Tang for
stimulating discussions and valuable advice at various stages of this
project.

\section{Notation I}
\label{notations}

\refmakepoint\textbf{Galois groups.}\label{notation-galois}
Let $F$ be a local field, $\kappa_F$ its residue field, and
$\F/F$ a Galois extension that splits the torus $T$ defined over $F$. Let
$\Gamma = \Gal(\F/F)$ and let $W_F$ and $W_{\F}$ be the Weil groups of $F$ and
$\F$. Recall that the relative Weil group of $\F/F$ is defined as $\W =
W_F/[W_{\F}, W_{\F}]$.

By the local class field theory, the relative Weil group is the
unique group extension
\[
	1\to \Ex \to \W \to \Gamma \to1
\]
corresponding to the fundamental class $\alpha\in H^2(\Gamma, \Ex)$.
Fix a system of representatives $\{w_\tau\}_{\tau\in\Gamma}$ for the
(right) cosets of $\Ex$. For any $g, h\in \W$, there is
a unique $\delta(g, h)\in \Ex$ and a unique $\tau\in\Gamma$ such that
\[
	g h = \delta(g, h) w_\tau
.\]
The assignment $(\tau, \sigma)\mapsto \delta(w_\tau, w_\sigma)$
is a cocycle which represents the fundamental class $\alpha$.

Let $p:\W\to\Gamma$ be the projection. We sometimes extend the
notation $w_\tau$ to allow subscripts be elements of $\W$, for
example, $w_g = w_{p(g)}$, and $w_{g\tau} = w_{p(g)\tau}$, for all $g\in\W,
\tau\in\Gamma$.

Let $U^{(n)}$ be the $n$-th congruence subgroup of $\F^\times$, so
$U^{(0)}=\mathcal O^\times_{\F}$ and $U^{(n)}=1 + \m^n_{\F}, n\geq 1$.
We define $W^{(n)}=\W/U^{(n)}$ and let $p_n:W^\np\to\Gamma$ be the
projection.

\makepoint\textbf{Torus and dual torus.}
Let $L = X^*(T)$ and $\hat L = X_*(T)$ be the weight lattice and
coweight lattice of $T$. They are both $\Gamma$-modules and therefore
$\W$-modules. The dual torus $\hat T = L \otimes \Gm$ is defined over
$\Z$, and so is the L-group $\cT = \hat T \rtimes \Gamma$.

We write $U_n = \hat L\otimes U^{(n)}$. They form a basis of
$\Gamma$-invariant open compact subgroups of $T(\F)$.

\makepoint[Picard stacks.]
All stacks in this paper will be fpqc-stacks. Let us review the
essentials of Picard stacks. The reader can refer to
\cite[\S~1.4]{SGA4Expose18} and \cite[Appendix~A]{CZ14}. Recall
that a Picard stack $\mathscr P$ over $S$ is a stack over $S$ together
with a bi-functor
\[
\otimes:\mathscr P\times\mathscr P\to\mathscr P,
\]
and the associativity and commutativity constraints
\[
a:\otimes\circ(\otimes\times 1)\cong\otimes\circ(1\times\otimes),
\quad
c:\otimes\cong\otimes\circ\mathrm{flip},
\]
such that for every $U$ over $S$, $\mathscr P(U)$ forms a Picard
groupoid (i.e. a symmetric monoidal groupoid in which every object
admits a monoidal inverse). A Picard stack is called strictly
commutative if $c_{x, x} = \id$ for every $x\in\mathscr{P}$. All
Picard stacks in this paper are strictly commutative.

Let $\Ch^\flat(S)$ be the category where objects are Picard stacks
over $S$ and morphisms are isomorphism classes of additive
functors.
In fact, $\Ch^\flat(S)$ is canonically enriched over itself. For
$\mathscr P_1, \mathscr P_2\in\Ch^\flat(S)$, we use
$\Homm(\mathscr P_1, \mathscr P_2)$ to denote the Picard stack
whose objects are additive functors from $\mathscr P_1$ to $\mathscr
P_2$ and morphisms are morphisms between additive functors
(cf. \cite[\S~1.4.7]{SGA4Expose18}). The dual of a Picard stack
$\mathscr P$ is defined to be $\mathscr P^\vee = \Homm(\mathscr P,
\BGm)$.

Let $C^{[-1, 0]}(S)$ be the category of 2-term complexes of 
sheaves of abelian groups $d: K^{-1}\to K^0$ over $S$.
Let $K\in C^{[-1, 0]}(S)$. We associate to it a Picard prestack
$\mathrm{pch}(K)$ whose section over $U$ is the following Picard
groupoid
\begin{itemize}
    \item Objects of $\mathrm{pch}(K)(U)$ are $K^0(U)$.
    \item If $x, y\in K^0(U)$, a morphism from $x$ to $y$ is an
    element $f\in K^{-1}(U)$ such that $df = y - x$.
\end{itemize}
Let $\ch(K)$ be the stackification of $\mathrm{pch}(K)$.

Let $D^{[-1, 0]}(S)$ denote the subcategory of the derived
category of abelian sheaves over $S$ consists of complexes $K$
with cohomology concentrated in degrees $-1$ and $0$.
Deligne proved that the functor $\ch$ induces an equivalence of
categories (cf. \cite[\S~1.4.15]{SGA4Expose18}):
\[
\ch: D^{[-1, 0]}\to\Ch^\flat(S).
\]
We choose a quasi-inverse of $\ch$, denoted by $\flat$, that sends
a Picard stack $\mathscr{P}$ to a complex
$\mathscr{P}^\flat \in D^{[-1, 0]}(S)$.

\makepoint[Representation functors.]
Let $S$ be a base scheme. Let $M$ be an affine group scheme over $S$
and let $G$ be an abstract group. We define $\repfun_{G, M}$ to be
the functor that sends a scheme $U$ over $S$ to $\Hom(G, M(U))$.

The functor $\repfun_{G, M}$ is represented by a scheme affine over the
base. This is because, if $I\subset G$ is a set of generators,
$\repfun_{G, M}$ is realized as a closed subscheme of $M^I$ defined by
the relations between the generators.

Let $G$ be a locally profinite group, and let $\{K_i\}_i$ be a basis
of open compact subgroups of $G$ such that $K_{i + 1}\leq K_i$. We
define
\[
\mathcal R^c_{G, M} = \varinjlim_i \mathcal R_{G / K_i, M}.
\]

The functor $\repfunc_{G, M}$ is represented by an ind-scheme because
each map $\repfun_{G/K_i, M}\to\repfun_{G/K_{i + 1}, M}$ is a closed
immersion.

\refmakepoint[Abelian sheaves.]
\label{abelian-sheaves}
All sheaves in this paper are fpqc-sheaves. For an abelian group $M$,
we write $\underline{M}$ for the constant abelian sheaf associated to
$M$.

For a locally profinite abelian group $N$, let $\{K_i\}_i$ be a
basis of open compact subgroups of $N$. We denote by
$N^\circ$ the sheaf
\[
        N^\circ = \varprojlim_i \underline{N/K_i}.
\]
Note that $N^\circ$ is \emph{not} a constant sheaf, because
the fpqc-topology remembers the locally profinite topology of $N$.
Indeed, for a test scheme $X = \Spec R$, we have
\[
N^\circ(X) = \cHom(\pi_0(X), N).
\]

If $N$ is profinite, each $\underline{N/K_i}$ is as a constant affine
group scheme, so $N^\circ$ is an affine group scheme (not necessarily
of finite type over the base). More generally, if $N$ is locally
profinite, $N^\circ$ is an ind-affine group scheme.

Let $T(F)_0$ denote the maximal compact subgroup of $T(F)$. For
convenience, we introduce the shorthand
\[
\uto = \big(T(F)_0\big)^\circ,\qquad \utf = T(F)^\circ.
\]
We write $\BTO$ and $\BTF$ for the classifying stack of $\uto$-torsors
and $\utf$-torsors.

\makepoint\textbf{Group cohomology.}
Let $G$ be a locally profinite group, and let $M$ be a discrete
abelian group with a continuous $G$-action.

Recall that $Z^1(G, M)$ is the abelian group of continuous 1-cocycles.
For every $m \in M$, the cocycle $\phi(g) = g m - m$ is continuous
becuase the action of $G$ on $M$ is continuous. This defines the
differential $d: M\to Z^1(G, M)$.
The first group cohomology is defined by $H^1(G, M) = Z^1(G, M)/d(M)$.

On the other hand, we define $\bar Z^1(G, M)$ to be the abelian group
of 1-cocycles when $G$ is viewed as an abstract group, and we define
$\bar H^1(G, M) = \bar Z^1(G, M)/d(M).$

Let $S = \Spec{R}$ be a base scheme, and let $\mathscr M$ be a
commutative group scheme over $S$ with a continuous $G$-action. We
define $\hZ^1(G, \mathscr M)$ to be the presheaf that sends a
commutative algebra $R'$ over $R$ to $Z^1(G, \mathscr M(R'))$.
$\hZ^1(G, \mathscr M)$ is in fact a sheaf. Let $d: \mathscr M\to
\hZ^1(G, \mathscr M)$ be the homomorphism of abelian sheaves such that
for every commutative algebra $R'$ over $R$, $d(R')$ coincides with
the differential defined before $d: \mathscr M(R')\to Z^1(G,
\mathscr M(R'))$. We define $\hH^1(G, \mathscr M)$ to be the Picard stack
\[
\hH^1(G, \mathscr M) = \ch\left(
\cdots\to 0\to\mathscr M\xrightarrow{d} Z^1(G, \mathscr M)\to 0\to\cdots
\right).
\]

\refmakepoint \textbf{The stack of Langlands parameters.}
\label{def-loc} We define the ind-scheme of framed Langlands
parameters over $\Spec\Z$ by $\Locc^\square = \hZ^1(\W, \hat T)$.
To see $\Locc^\square$ is an ind-scheme, we can rewrite it as:
\[
	\Locc^\square = \varinjlim_n \Locc^{\square, (n)},\quad
	\text{where}\quad
	\Locc^{\square, (n)} = \repfun_{W^{(n)},\cT}
	\times_{\repfun_{W^{(n)}, \Gamma}}
	\{p_n\}
.\]

The stack of Langlands parameters is defined to be the quotient stack
$\Locc = [\Locc^\square/\hat T]$. Because the action of $\hat T$
stabilizes each $\Locc^{\square, (n)}$, $\Locc$ is an ind-algebraic
stack:
\[
        \Locc = \varinjlim_n \Locc^{(n)}, \quad\text{where}\quad
        \Locc^{(n)} = \Locc^{\square, (n)}/\hat T.
\]
It is easy to see $\Locc = \hH^1(\W, \hat T)$. This endows
a Picard stack structure on $\Locc$, as explained above.

\refmakepoint \textbf{The stack} $\Tor$ \textbf{.}
\label{def-tor} We follow the setup in \cite{Kot14}. Let
\begin{equation*}
        \xi: \hat L \to \Hom(\F^\times, \F^\times\otimes\hat L)
	= Z^1(\Ex, T(\F))
\end{equation*}
be the adjoint of the identity. Taking $\Gamma$-invariants we get
\[
        \xi^\Gamma: \hat L^\Gamma\to Z^1(\Ex, T(\F))^\Gamma
\]
We form the set of algebraic cycles by the following fibre product:
\[
\begin{tikzcd}
Z^1_{\mathrm{alg}}(\W, T(\F)) \ar[r]\ar[d]& \hat L^\Gamma \ar[d]\\
Z^1(\W, T(\F)) \ar[r, "\mathrm{res}"]
& Z^1(\Ex, T(\F))^\Gamma.
\end{tikzcd}
\]
In other words, the algebraic cycles are those 1-cycles whose
restriction to $\F^\times\subset \W$ is an algebraic character of $T$
(coming from the restriction, this character is automatically
$\Gamma$-invariant). The differential $d: T(E)\to Z^1(\W, T(E))$
factors through $d': T(E)\to Z^1_\alg(\W, T(E))$, because the
restriction to $\Ex\subset\W$ is always the trivial character.

In \textit{loc.\,cit.}, Kottwitz proved that the cokernel of $d'$ is
isomorphic to the Kottwitz set $B(T) = \hat L_\Gamma$. Thus we have a
long exact sequence
\[
0\to T(F) \to T(E) \xrightarrow{d'} Z^1_\alg(\W, T(E)) \to \hat L_\Gamma\to 0.
\]
Since $T(F)$ is a closed subgroup of $T(E)$, the image of $d'$ is
endowed with the quotient topology which is locally profinite.
We endow $Z^1_\alg(\W, T(E))$ with the unique topology such that it
coincides with the quotient topology on $\Im(d')$, and the projection
to $\hat L_\Gamma$ is continuous (the latter is equipped with
the discrete topology). This makes $Z^1_\alg(\W, T(E))$ a locally
profinite group, and the map $d'$ is continuous.

Now we can define $\Tor$ as:
\begin{equation*}
        \Tor = \mathrm{ch}\left(
	\cdots\to0\to
	T(\F)^\circ\to
	Z^1_\alg(\W, T(E))^\circ\to
	0\to\cdots\right).
\end{equation*}

\begin{remark}
\label{decomposition-tor}
\begin{enumerate}[label=(\roman*)]
\item 
The isomorphism classes of the groupoid $\Tor(\Z)$ is the Kottwitz set
$B(T) = \hat L_\Gamma$. The automorphism group of the identity object
is $T(F)$.
\item
It is easy to see that
\[
\Tor = \bigsqcup_{\hat L_\Gamma} \, [*/T(F)^\circ].
\]
Although the definition of $\Tor$ involves inverse limit, we mostly
consider it as an ind-stack by considering $\Tor$ as the union of its
connected components.
\item
Let $\breve F$ be the completion of a maximal unramified
extension of $F$, and let $\sigma$ be a Frobenius element. $\Tor(\Z)$
is isomorphic to the groupoid of pairs $(\mathcal E, \phi)$ where
$\mathcal E$ is a $T$-torsor over $\breve F$ and $\phi:\mathcal
E\cong\sigma^*\,\mathcal E$ is an isomorphism of $T$-torsors. By a
Tannakian formalism, this is the groupoid of exact tensor functors
from $\mathrm{Rep}_T$ to the monoidal category of isocrystals over
$F$. This justifies the notation $\Tor$.
\end{enumerate}
\end{remark}

\section{Lemmas in Homological Algebra}
\label{lemmas-in-HA}
\noindent
In this section and the next section, we establish some lemmas
in an abstract setting.

Let $\Gamma$ be a finite group, and let $A$ be an abelian group
equipped with a $\Gamma$-action. Let $\alpha\in H^2(\Gamma, A)$
be the cohomology class representing a group extension
\[
1\to A\to G\to \Gamma\to 1.
\]
\begin{definition}
\label{tate-nakayama-criterion}
We say the triple $(\Gamma, A, \alpha)$ satisfies the Tate-Nakayama
criterion if for every subgroup $\Gamma'\leq\Gamma$, and for all
$r\in\Z$, taking the cup product with
$\Res_{\Gamma/\Gamma'}(\alpha)$ induces isomorphisms
\[
\hat H^r(\Gamma', \Z) \to \hat H^{r + 2}(\Gamma', A).
\]
We remark that this notion equivalent to $(\Z, A, \alpha)$ being a
Tate-Nakayama triple in the sense of Kottwitz \cite{Kot14}.
\end{definition}

\begin{prop}
\label{triples}
Let $E/F$ be a finite extension of non-Archimedean local fields.
The following triples satisfy the Tate-Nakayama criterion.
\begin{enumerate}[label=(\roman*)]
        \item
        $\Gamma = \Gal(E/F)$, $A = \Ex$ and the fundamental class
        from class field theory $\alpha\in H^2(\Gamma, \Ex)$.
        \item
        Assume $E/F$ is tamely ramified. Let $\Gamma = \Gal(E/F)$,
        $A = \Ex/U^{(n)}$ and let $\alpha$ be the image of the
        fundamental class under
        \[
                H^2(\Gamma, \Ex)\to H^2(\Gamma, \Ex/U^{(n)}).
        \]
        \item 
        Let $\Gamma = I$ be the inertia group of $E/F$, let $A =
        I^\ab_E$ be the abelianization of the absolute inertia group of $E$,
        and let $\alpha = u_{\breve E/\breve F}\in H^2(I, I^\ab_E)$
        be the canonical class (for notation, see Section~\ref{inertia-group}).
\end{enumerate}
\end{prop}
\begin{proof}
Part (i) of the proposition is standard.

We move on to part (ii). Notice that the congruence subgroups
$U^{(n)}$ are isomorphic to the additive group $\mathfrak{m}_E^n$ via
logarithm and in turn are isomorphic to $\mathcal O_E$. When $E/F$ is
tamely ramified, it is well known that $O_E$ is free of rank $1$ as an
$O_F[\Gamma]$-module, hence it is cohomologically trivial. As a
result, for every subgroup $\Gamma'\subset \Gamma$, the natural map
$\hat H^r(\Gamma', \Ex)\to\hat H^r(\Gamma', \Ex/U^\np)$ is an
isomorphism for $r\geq 1$. By standard argument, this implies
isomorphism for all $r\in\Z$.

Part (iii) is the content of \cite[\S~2.5, Proposition~9]{Ser61}. The reader
can also refer to \cite[\S~2.5]{DW23}.
\end{proof}

The following lemma appears to be new.

\begin{lem}
Let $A, G, \Gamma, \alpha$ be as above. Let $M$ be a $\Gamma$-module,
or equivalently, a $G$-module on which $A$ acts trivially. Then we
have a commutative diagram
\[
\hspace{-8ex}
\begin{tikzcd}
& H_2(\Gamma, M) \ar[r, "d_2"] \ar[d, "\cup\alpha"]
& H_1(A, M)_\Gamma \ar[r] \ar[d, equal]
& H_1(G, M) \ar[r] \ar[d, "\mathrm{res}"]
& H_1(\Gamma, M) \ar[r] \ar[d, "\cup\alpha"] & 0\\
0 \ar[r] & \hat{H}^{-1}(\Gamma, M\otimes A) \ar[r] &
(M\otimes A)_\Gamma \ar[r] &
(M\otimes A)^\Gamma \ar[r] &
\hat{H}^0(\Gamma, M\otimes A) \ar[r] & 0.
\end{tikzcd}
\]
The first row of the diagram is the long exact sequence associated to the
Lyndon-Hochschild-Serre spectral sequence. The second row is the definition
of the Tate cohomology. The first and last vertical maps are the cup
product with $\alpha$. The map $\mathrm{res}$ is the restriction map
$H_1(G, M)\to H_1(A, M)^\Gamma$.

In particular, if $(\Gamma, A, \alpha)$ satisfies the Tate-Nakayama
criterion and $M$ is torsion-free as an abelian group, we have
\[
H_1(G, M) \cong (M \otimes A)^\Gamma.
\]
\end{lem}
\begin{proof}
    See Appendix A. \qedhere
\end{proof}

\begin{cor}
\label{lan-result}
\begin{enumerate}[label=(\roman*)]
    \item
    $H_1(\W, \hat L)\cong T(F)$. This is first proved by Langlands in
    \textup{loc.\,cit.}
    \item
    If $E/F$ is tamely ramified, then
    $H_1(\Wn, \hat L)\cong (T(E)/U_n)^\Gamma = T(F)/(U_n \cap T(F))$.
    \item
    $H_1(\I, \hat L)\cong (\hat L\otimes I^\ab_E)^I$.
\end{enumerate}
\end{cor}

Even when $E/F$ is not tamely ramified, the group $H_1(\Wn, \hat L)$
remains isomorphic to a quotient of $T(F)$ by a compact open subgroup.
We just lack of an explicit description of this subgroup.

\begin{prop}
    We have $T(F) = \varprojlim H_1(\Wn, \hat L)$, and each
    $H_1(\Wn, \hat L)$ has the discrete topology.
\end{prop}
\begin{proof}
    By the Lydon--Hochschild--Serre spectral sequence, we have
    \[
    H_1(U^\np, \hat L)_\Gamma\xrightarrow{\cores} H_1(\W, \hat L)\to
    H_1(\Wn, \hat L)\to 0,
    \]
    where the middle term $H_1(\W, \hat L)\cong T(F)$.
    Since $H_1(\Wn, \hat L)$ is finitely generated as an abstract
    abelian group, it must be a discrete quotient of $T(F)$.
    Since $H_1(U^\np, \hat L)\cong U_n\subset T(E)$, and
    the map $\cores: (U_n)_\Gamma \to T(F)$ is the norm map, it is
    clear that $\bigcap_n \mathrm{Im}(\cores) = 1$, therefore $T(F) =
    \varprojlim H_1(\Wn, \hat L)$.
\end{proof}

\begin{definition}
    \label{def-vn} 
    We define $V_n = \ker(H_1(\W, \hat L)\to H_1(\Wn, \hat L))$. By
    the previous proposition, they are a basis of open compact
    subgroups of $T(F)$. In \parti, if $E/F$ is tamely ramified, $V_n
    = U_n \cap T(F)$ by Corollary \ref{lan-result}.
\end{definition}

The next lemma is very general but important to the proof
of Lemma \ref{cont-iff}.
\begin{lemma}
\label{open-map}
Let $X, Y$ be locally profinite abelian groups, and let $f:X\to Y$ be
any group homomorphism. Assume that for one (hence every) compact open
subgroup $U\subset X$, the quotient $X/U$ is a finitely generated
abelian group. Assume that $f$ has finite cokernel. Then $f$ is open.
\end{lemma}
\begin{proof}
Let $U$ be any compact open subgroup of $X$, we show that $f(U)$ is
open. We have a commutative diagram with exact rows and columns:
\[
\begin{tikzcd}
    1 \ar[d] & 1 \ar[d] &&\\
    U \ar[r] \ar[d] & f(U) \ar[d] &&\\
    X \ar[r]\ar[d] & Y \ar[r]\ar[d]
    & Z \ar[r]\ar[d, twoheadrightarrow, "g"] & 1\\
    X' \ar[r]\ar[d] & Y' \ar[r]\ar[d]
    & Z' \ar[r] & 1\\
    1 & 1 &&\\
\end{tikzcd}
\]
The map $g$ is defined by lifting an element $z\in Z$ to $y\in Y$, and
sends it through $Y\to Y'\to Z'$. It is straightforward to check that
$g$ is well-defined and that $g$ is a surjection.

Because $Z$ is finite, so is $Z'$. Since $X'$ is a finitely generated
abelian group, so is $Y'$. Therefore $f(U)$ is open.
\end{proof}

\makepoint From now on, we assume that $A$ and $G$ are locally
profinite. Let $L$ be a free abelian group of finite rank with a
$\Gamma$-action, and let $\hat L$ be its dual.

The homology groups
$H_1(A, \hat L)$ and $H_1(G, \hat L)$ have natural topologies described
below. $H_1(A, \hat L)\cong A\otimes \hat L$ has the product topology.
We have the Lyndon--Hochschild--Serre spectral sequence for group
homology:
\[
H_2(\Gamma, \hat L)\to
H_1(A, \hat L)_\Gamma\xrightarrow{\cores'} H_1(G, \hat L)
\to H_1(\Gamma, \hat L)\to 0.
\]
$H_1(A, \hat L)_\Gamma$ is the quotient of $H_1(A, \hat L)$ by a
closed subgroup, hence with the quotient topology it is locally
profinite. Now since $H_i(\Gamma, \hat L)$ are finite for $i = 1, 2$,
the image of $H_1(A, \hat L)_\Gamma$ under $\cores'$ with the quotient
topology is locally profinite. This uniquely determines a locally
profinite topology on $H_1(G, \hat L)$.

The corestriction map
\[
        \cores: H_1(A, \hat L)\to H_1(G, \hat L)
\]
is continuous by definition. It has finite cokernel, hence it is open
by Lemma \ref{open-map}.

\begin{lem}
For any commutative ring $R$, we have a short exact sequence 
\begin{equation}
	\label{coarse}
	0 \to \Ext^1(\hat L_\Gamma, \Rx) \xrightarrow{i}
	\bar H^1(G, L\otimes \Rx) \xrightarrow{p}
	\Hom(H_1(G, \hat L), \Rx) \to 0
.\end{equation}
\end{lem}
\begin{proof}
Let $C_\bullet \to \Z \to 0$ be the bar
resolution of the trivial module $\Z$ by free $\Z[G]$-modules.
The group cohomology $\bar H^1(G, L\otimes\Rx)$ is exactly the first
cohomology of the following complex: 
\begin{equation*}
	\Hom_G(C_\bullet, L\otimes R^\times)
	= \Hom((C_\bullet \otimes \hat L)_G, R^\times)
.\end{equation*}
Because the complex $(C_\bullet \otimes \hat L)_G$ is a complex of
free abelian groups, and its homology calculates the group homology of
the $G$-module $\hat L$, we apply the universal coefficient theorem
and get
\[
	0\to\Ext^1(H_0(G, \hat L), R^\times)
	\to \bar H^1(G, L\otimes\Rx)
	\to \Hom(H_1(G, \hat L), R^\times)
	\to 0
.\]
\end{proof}

\begin{lemma}
We have a short exact sequence
\[
0 \to \Ext^1(\hat L_\Gamma, R^\times) \xrightarrow{i'}
H^1(G, L\otimes\Rx) \xrightarrow{p'}
\cHom(H_1(G, \hat L), R^\times) \to 0
.\]
\label{coarse-cont}
\end{lemma}
\begin{proof}
We first show the image of $\Ext^1(\hat L_\Gamma, R^\times)$ is
contained in $H^1(G, L\otimes\Rx)$.
To do this, we work out an explicit formula for the map $i'$.
Let $\hat L_0 = \{gx - x\,|\,g\in \Gamma, x\in \hat L\}$, so we have
\[
	0 \to \hat L_0 \to \hat L \to \hat L_\Gamma \to 0
.\]
This short exact sequence gives rise to
\[
	\Hom(\hat L_0, R^\times) \xrightarrow{\alpha}
	\Ext^1(\hat L_\Gamma, R^\times)
	\to 0
.\]
On the other hand, we have a map $\Hom(\hat L_0, R^\times) \to
\Hom(G\times\hat L, \Rx)$ by sending $\phi\in\Hom(\hat L_0, R^\times)$ to
\[
	(g, x) \mapsto \phi(g^{-1}x - x)
.\]
By adjunction this defines a cocycle $G \to L \otimes
R^\times$. Denote this map by $\beta: \Hom(\hat L_0, R^\times) \to
\bar H^1(G, L\otimes\Rx)$. Then $i' = \beta \circ \alpha^{-1}$,
independent of the choice of a preimage of $\alpha$. The image of $i'$
clearly consists of continuous cocycles because they factor through
$G\to\Gamma$.

It remains to show that for $\phi\in\bar H^1(G, L\otimes\Rx)$, it is
a continuous cocycle if and only if $p(\phi)$ is a continuous
homomorphism. The map $p$ is induced by $\pi$ on the level of cocycle
as follows:
\[
\begin{tikzcd}
	\bar Z^1(G, L\otimes\Rx) \ar[r, "\pi"] &
	\Hom(H_1(G, \hat L), R^\times)\\
	\phi:G \to L\otimes R^\times \ar[r, mapsto] &
	\Big(
	\sum g\otimes x \mapsto
	\sum \langle \phi(g), x \rangle
	\Big),
\end{tikzcd}
\]
where $\phi$ is a cocycle, $\sum g\otimes x \in G\otimes \hat L$ is a
cycle representing an element in $H_1(G, \hat L)$, and the pairing
$\langle -,- \rangle$ is the evaluation of $\hat L$ on $L$. It remains
to show that $\phi$ is continuous if and only if $\pi(\phi)$ is. This
is the content of the next lemma.
\end{proof}

\begin{lemma}
\label{cont-iff}
$\phi$ is continuous if and only if $\pi(\phi)$ is continuous.
\end{lemma}
\begin{proof}
On the one hand, since $\cores:H_1(A, \hat L) \to H_1(G, \hat L)$ is
both continuous and open, $\psi:H_1(G, \hat L)\to R^\times$ is
continuous if and only if $\psi\circ\cores$ is.

On the other hand, it is clear that a cocycle is continuous
if and only if its restriction to $A$ is. We have the following
commutative diagram:
\[
\begin{tikzcd}
	\bar Z^1(G, L\otimes\Rx) \ar[r, "\pi"]\ar[d]&
	\Hom(H_1(G, \hat L), R^\times) \ar[d]\\
	\bar Z^1(A, L\otimes\Rx) \ar[r, "\tau"]&
	\Hom(H_1(A, \hat L), R^\times),
\end{tikzcd}
\]
where the map $\tau$ sends a cocycle $f$ to an homomorphism
$\tau(f):H_1(A, \hat L)\to R^\times$, $a\otimes x\mapsto\langle
f(a), x\rangle$. It is then clear that $f$ is continuous if
and only if $\tau(f)$ is, and the lemma follows.
\end{proof}

In the rest of this section, we introduce a lemma of Kottwitz
\cite[\S~3]{Kot14}. First, we need some notation.

\makepoint \textbf{Definition of \textmd{$Z^1_Y(G, M)$} and
\textmd{$H^1_Y(G, M)$}.} Let $A, G, \Gamma, \alpha$ be as above. We
consider a triple $(M, Y, \xi)$ where $M$, $Y$ are $\Gamma$-modules
and $\xi:Y\to\Hom(A, M)$ is a $\Gamma$-map. Taking $\Gamma$-invariants
gives us a map $\xi^\Gamma: Y^\Gamma\to\Hom_\Gamma(A, M)$. We define
$Z^1_Y(G, M)$ as the fibre product:
\[
\begin{tikzcd}
    Z^1_Y(G, M) \ar[r]\ar[d] & Y^\Gamma \ar[d, "\xi^\Gamma"]\\
    Z^1(G, M) \ar[r, "\res"] & \Hom_\Gamma(A, M),
\end{tikzcd}
\]
and we define $H^1_Y(G, M)$ as the fibre product:
\[
\begin{tikzcd}
    H^1_Y(G, M) \ar[r]\ar[d, "i"] & Y^\Gamma \ar[d, "\xi^\Gamma"]\\
    H^1(G, M) \ar[r, "\res"] & \Hom_\Gamma(A, M).
\end{tikzcd}
\]

Let $B^1_Y(G, M)$ stands for the subgroup of $Z^1_Y(G, M)$ consisting
of pairs $(0, m)$ where $m$ is a 1-coboundary of the $G$-module $M$;
$(0, m)$ defines an element in the fibre product because $A$ acts
trivially on $M$, therefore the restriction of $m$ to $A$ is trivial.
There is a natural map $Z^1_Y(G, M)\to H^1_Y(G, M)$ and it induces
an isomorphisms $H^1_Y(G, M) = Z^1_Y(G, M)/B^1_Y(G, M)$.

\makepoint[A lemma of Kottwitz.] 
We define a map $c': Y\to Z^1_Y(G, M)$ as follows. It suffice
to define two maps $Y\to Y^\Gamma$ and $Y\to Z^1(G, M)$ and
check the commutativity in the product diagram.
Let the first map be the norm map $N:Y\to Y^\Gamma$, and the
second to be the composition
\[
Y\xrightarrow{\xi} Z^1(A, M)
\xrightarrow{\inf} Z^1(G, M).
\]
To check that $\xi^\Gamma\circ N = \res\circ\inf\circ\xi$, it suffice
to note that $\res\circ\inf$ coincides with the norm map
$N: Z^1(A, M)\to Z^1(A, M)$.
Postcomposing $c'$ with the projection $Z^1_Y(G, M)\to H^1_Y(G, M)$
yields a map $c_0: Y\to H^1_Y(G, M)$.
By Lemma 3.3 in \cite{Kot14}, $c_0$ factors through
$c: Y_\Gamma\to H^1_Y(G, M)$.

There is a tautological pairing $A\otimes\Hom(A, M)\to M$, which,
combined with the given map $\xi:Y\to\Hom(A, M)$ gives a pairing
$A\otimes Y\to M$. The cup product with $\alpha\in H^2(\Gamma, A)$
gives
\begin{equation}
    \hat H^i(\Gamma, Y)\to\hat H^{i + 2}(\Gamma, M).
    \label{cup-product}
\end{equation}

The following is Lemma 3.5 from \cite{Kot14}.

\begin{lem}
\label{lem-Kottwitz}
\begin{enumerate}[label=(\roman*)]
  \item The following diagram commutes:
  \[
  \begin{tikzcd}
    H^1(\Gamma, M) \arrow[r, "i"] & H^1_Y(G, M) \\
    H^{-1}(\Gamma, Y) \arrow[u, "\cup\alpha"] \arrow[r] & Y_\Gamma. \arrow[u, "c"]
  \end{tikzcd}
  \]

  \item The homomorphism $c : Y_\Gamma \to H^1_Y(G, M)$ is an isomorphism
  if and only if the map in (\ref{cup-product}) is bijective for
  $i = -1$ and injective for $i = 0$.
\end{enumerate}
\end{lem}

Let $L$ be a free abelian group of finite rank, and let $\hat L$ be
its dual. We take $M = \hat L\otimes A$, $Y = \hat L$ and let
$\xi$ be the adjoint of the identity
\[
\xi: \hat L\to\Hom(A, \hat L\otimes A).
\]
In this case, we denote $Z^1_Y(G, M)$ and $H^1_Y(G, M)$ by
$\Zalg(G, M)$ and $\Halg(G, M)$. 

If $(\Gamma, A, \alpha)$ satisfies the Tate-Nakayama criterion,
then the maps in (\ref{cup-product}):
\[
\hat H^i(\Gamma, \hat L)\to\hat H^{i + 2}(\Gamma, \hat L\otimes A)
\]
are isomorphisms for all $i\in\Z$. Therefore we can apply Lemma
\ref{lem-Kottwitz} to obtain:

\begin{cor}
    \label{kot-result}
    The map $c:\hat L_\Gamma\to \Halg(G, \hat L\otimes A)$
    is an isomorphism.
\end{cor}

\section{Split Exact Sequence}
\label{short-exact-sequence}
\noindent
Let $A, G, \Gamma, \alpha$ be the same as the last section.
\begin{prop}
    \label{picard-exact-seq}
    We have a short exact sequence of Picard stacks
    \[
    1\to\Homm(\hat L_\Gamma, \BGm)\to
    \hH^1(G, L\otimes\Gm)\to
    \rfc{H_1(G, \hat L), \Gm}\to 1.
    \]
\end{prop}
\begin{proof}
Recall that the functor $\mathrm{ch}:D^{[-1, 0]}(S) \to \Ch(S)^\flat$
is an equivalence of categories. Let $(-)^\flat$ be a quasi-inverse. By
definition, we need to show there exists an exact triangle
\begin{equation}
    \label{exact-triangle}
	\Homm(\hat L_\Gamma, \BGm)^\flat \to
	\hH^1(G, L\otimes\Gm)^\flat \to
	(\rfc{H_1(G, \hat L), \Gm})^\flat \xrightarrow{+1}.
\end{equation}

We first define the homomorphisms
\[
	\Homm(\hat L_\Gamma, \BGm)^\flat \xrightarrow{\bar \beta}
	\hH^1(G, L\otimes\Gm)^\flat \xrightarrow{\bar \pi}
    (\rfc{H_1(G, \hat L), \Gm})^\flat
\]
as follows, where each column is a 2-term complex of abelian sheaves
corresponding to the Picard stack above:
\[
\hspace{-15ex}
\begin{tikzcd}
	{\text{deg}=0} & \Homm(\hat L_0, \Gm) \ar[r, "\beta'"]&
	\hZ^1(G, L\otimes\Gm) \ar[r, "\pi'"] &
	\rfc{H_1(G, \hat L), \Gm} \\
	{\text{deg}=-1} & \Homm(\hat L, \Gm) \ar[u, "\mathrm{res}"]
	\ar[r, "\sim"] &
	L\otimes\Gm \ar[u, "d"] \ar[r] &
	0. \ar[u]
\end{tikzcd}
\]
For any commutative ring $R$, we define the map $\beta'(R)$ to be
$\beta$ as defined in Lemma \ref{coarse-cont}. The map $\pi'(R)$
is the restriction of $\pi$ to the continuous cocycles. It is
immediate to check that the squares commute, so they define
homomorphisms of complexes.

We show that the three terms form an exact triangle by showing
a quasi-isomorphism $\cofib(\bar\beta)\cong\rfc{H_1(G, \hat L), \Gm}$.
Note $\mathrm{cofib}(\bar\beta)$ is nothing other than the following
complex:
\[
	\Homm(\hat L, \Gm) \xrightarrow{(\id, -\mathrm{res})}
        \left(L\otimes\Gm\right) \oplus \Homm(\hat L_0, \Gm)
	\xrightarrow{\left(\begin{smallmatrix} d\\
	\beta'\end{smallmatrix} \right)}
	\hZ^1(G, L\otimes\Gm)
.\]
Because the first map is injective, $H^{-2}(\mathrm{cofib}(\bar\beta))
= 0$. Because $d \circ \id = \beta'\circ\mathrm{res}$ and
$\beta'$ is injective, $H^{-1}(\mathrm{cofib}(\bar\beta)) = 0$.

We know that \(
H^0(\mathrm{cofib}(\bar\beta)) = Z^1(G, L\otimes\Gm) /
\Im \left(\begin{smallmatrix} d\\ \beta'\end{smallmatrix}\right)
\).
To evaluate this, we first take the quotient \(
Z^1(G, L\otimes\Gm) /
\Im' \left(\begin{smallmatrix} d\\ \beta'\end{smallmatrix}\right)
\)
as presheaf, where $\Im'$ is the image of maps of presheaves. We have
\begin{align*}
	&\quad\, Z^1(G, L\otimes\Rx)/(\Im'(d) + \Im'(\beta')) \\
        &= H^1(G, L\otimes\Rx)
	/(\Im'(\beta')/(\Im'(\beta')\cap \Im'(d)))\\
	&= H^1(G, L\otimes\Rx)/\Ext^1(\hat L_\Gamma, R^\times)\\
	&= \cHom(H_1(G, \hat L), R^\times)
.\end{align*}
\end{proof}

\begin{prop}
\label{splitting}
The +1 map of the exact triangle (\ref{exact-triangle}) is zero.
\end{prop}
\begin{proof}
Recall that the +1 map is the inverse of the quasi-isomorphism
showed above composed with a canonical projection
\[
	\rfc{H_1(G, \hat L), \Gm}
	\overset{\sim}{\leftarrow}\cofib(\bar\beta)
	\to \Homm(\hat L_\Gamma, \BGm)[1]
.\]

By the universal coefficient theorem, there exists a
collection of maps
\begin{equation*}
	\rho: \Hom(H_1(G,\hat L), R^\times) \to
	\bar Z^1(G, L\otimes\Rx)
\end{equation*}
that is functorial in $R$, such that after composing with the
projection $\bar Z^1(G, L\otimes\Rx) \to \bar H^1(G, L\otimes\Rx)$,
give splittings of the following surjections:
\[
	\pi: \bar H^1(G, L\otimes\Rx) \to
    \Hom(H_1(G, \hat L), R^\times)
.\]
For a continuous homomorphism $f:H_1(G, \hat L) \to R^\times$,
$\rho(f)$ is automatically a continuous cocycle because its orbit
under action of $L\otimes\Rx$ must contain one continuous cocycle, but
since every coboundry is continuous, $\rho(f)$ is continuous.
Therefore we have a family of maps that is functorial in $R$:
\begin{equation*}
	\cHom(H_1(G, \hat L), R^\times) \to
	Z^1(G, L\otimes\Rx)
.\end{equation*}

These maps give rise to an inverse of the quasi-isomorphism
$\mathrm{cofib}(\bar\beta) \to \rfc{H_1(G, \hat L),\Gm}$ by the
composition
\[
	\rfc{H_1(G, \hat L), \Gm} \to \hZ^1(G, L\otimes\Gm)
	\to \mathrm{cofib}(\bar\beta)
.\]
Composing this inverse with $\mathrm{cofib}(\bar\beta)\to
\Homm(\hat L_\Gamma, \BGm)[1]$ is nothing but the +1 map, but it is
also zero by construction.
\end{proof}

\begin{cor}
\label{two-ses}
For an arbitrary torus $T$ over a non-Archimedean local field $F$, we
have the following split short exact sequences of Picard stacks
\begin{gather}
	1 \to \Homm(\hat L_\Gamma, \BGm) \to
	\Locc \to \rfc{T(F), \Gm} \to 1,\label{ses}\\
    1 \to \Homm(\hat L_\Gamma, \BGm) \to
    \Locn \to \repfun_{T(F)/V_n, \Gm} \to 1.\label{ses2}
\end{gather}
\end{cor}
\begin{proof}
    Let $L$ be the character lattice of $T$ and let $\hat L$ be its dual. In
    Proposition \ref{picard-exact-seq} and Proposition \ref{splitting},
    substituting $G$ for $\W$, $\hH^1(G, L\otimes\Gm)$ for $\Locc$ and $H_1(G,
    \hat L)$ for $T(F)$, or substituting $G$ for $\Wn$, $\hH^1(G, L\otimes\Gm)$
    for $\Locn$ and $H_1(G, \hat L)$ for $T(F)/V_n$ (for definition
    of $V_n$, see Definition \ref{def-vn}).
\end{proof}

Let $S$ be an arbitrary base scheme, and let $M$ be a finitely generated
abelian group over $S$. That is to say, locally on $S$, $M$ is isomorphic
to a constant sheaf of a finitely generated abelian group. Recall that
the Cartier dual of $M$, denoted by $D(M)$, is the sheaf that associates
to a scheme $U$ over $S$ the set of group scheme homomorphisms
$\Hom(M\times_S U, \Gm)$.

\begin{rmk}
(i)  The first term in both short exact sequences has the alternative
form 
\[
        \Homm(\hat L_\Gamma, \BGm) = \B\hat T^\Gamma.
\]
This is because $X^*(\hat T^\Gamma) = \hat L_\Gamma$. By Example
A.3.3. in \cite{CZ14}, we have
\[
        \Homm(\hat L_\Gamma, \BGm) \cong \B(D(\hat L_\Gamma))
        \cong \B\hat T^\Gamma.
\]

(ii) By the previous remark, the sheaf of automorphism of the identity
object in $\Locc$ is the Picard dual of the sheaf of isomorphism
classes of $\Tor$, namely
\[
(\B\hat T^\Gamma)^\vee\cong\hat L_\Gamma.
\]
On the other hand, the sheaf of isomorphism classes of $\Locc$ is the
Picard dual of the automorphism group of the identity object in
$\Tor$, namely 
\[
(\rfc{T(F), \Gm})^\vee\cong\BTF.
\]
This is a prelude to the duality
between $\Locc$ and $\Tor$, which is the theme of the next section.
\end{rmk}

\section{Duality}
\label{duality}
\noindent
The first main result of this paper is the following.
\begin{thm}
	\label{main}
	There is a unique family of Poincar\'e line bundles
	\[
		\mathcal L_T:	\Tor\times\Locc \to \BGm
	\]
	for every torus $T$ over $F$ that satisfies the three conditions
	below. When $2$ is invertible in $S$, these Poincar\'e line bundles
	induce isomorphisms $\Tor\cong\Locc^\vee$ for every torus $T$.

\begin{enumerate}[label=(\alph*)]

	\item \textbf{Functorality.} Let $g: S\to T$ be a map between torus.
	It induces $\alpha:\Locc \to \Loc_{^c S, F}$ and $\beta:\TorS \to
	\Tor$. The following two line bundles on $\TorS\times\Locc$ are
	canonically isomorphic:
	\[
		(\beta\times\id)^*\, \mathcal L_T = (\id\times\alpha)^*\, \mathcal L_S.
	\]

	\item \textbf{Split case.} For the split torus $T = \Gm$, both $\Tor$
	and $\Locc$ canonically split as
	\begin{align*}
		\Tor &= \Z\times \B F^\times,\\
		\Locc &= \rfc{F^\times, \Gm}\times \BGm.
	\end{align*}
	Then $\mathcal L_T$ is the canonical line bundle $\Tor\times\Locc\to
	\BGm$ via the tautological pairings
	\begin{align*}
		\rfc{F^\times, \Gm}\times
		\B F^\times &\to \BGm\\
		\Z\times \BGm&\to \BGm.
	\end{align*}

	\item \textbf{Induced case.} For an induced torus $S =
	\Res_{E/F}\Gm$, let $T$ be the split torus over $E$, the Shapiro
	isomorphism $\TorS\times\Loc_{{}^L S, F} \cong \TorTE\times
    \Loc_{\cT, E}$ identifies the line bundles
    $\mathcal L_S$ and $\mathcal L_T$.

\end{enumerate}

\end{thm}

Our strategy is to introduce an auxiliary stack $\T_{T/F}$ for every
torus $T$ over $F$ in place of $\Tor$. When it is clear from the
context, we will write $\T_T$ or even $\T$ to stand for $\T_{T/F}$. We
will first define a family of line bundles on $\T \times\Locc$. Then
we verify the main theorem for $\T$ instead of $\Tor$. Finally, we
will show a canonical isomorphism $\T\cong\Tor$.

\makepoint\textbf{Definition of $\T$.}
Let $C_\bullet(\W, \hat L)$ be the bar resolution that
calculates the group homology of the $\W$-module $\hat L$.
The boundaries of this chain complex are defined by
\[
	B_i(\W, \hat L) = \Im\left(C_{i+1}(\W, \hat L)
	\xrightarrow{d}
	C_i(\W, \hat L) \right).
\]
Consider the following homomorphism
\[
\tfrac{C_1(\W, \hat L)}{B_1(\W, \hat L)}
\xrightarrow{d} C_0(\W, \hat L).
\]
By definition, $C_0(\W, \hat L) = \hat L$ and the image of
$d$ is $\hat L_0 = \{gx - x\,|\,g\in \Gamma, x\in \hat L\}$.
Therefore, we have a short exact sequence
\begin{equation}
\label{hello}
0\to H_1(\W, \hat L)\to \tfrac{C_1(\W, \hat L)}{B_1(\W, \hat L)}
\to\hat L_0\to 0.
\end{equation}
This shows that we can endow $\frac{C_1(\W, \hat L)}{B_1(\W, \hat L)}$
with the unique topology making the inclusion
$H_1(\W, \hat L)\cong T(F)$ an open embedding.

Now we can define
\[
	\T :=
	\mathrm{ch}\left(
	\cdots\to 0\to
	\tfrac{C_1(\W, \hat L)}{B_1(\W, \hat L)}^\circ
	\xrightarrow{d} \underline{\hat L}
	\to0\to\cdots
	\right)
.\]
The groupoid $\T(\Z)$ has isomorphism classes indexed by $\hat L_\Gamma$
and the automorphism group of the identity is $H_1(\W, \hat L)\cong T(F)$.

\begin{prop}
    There is a split short exact sequence
    \begin{equation}
    \label{tor-ses}
    1\to \BTF\to\T\to\underline{\hat L_\Gamma}\to 1.
    \end{equation}
\end{prop}
\begin{proof}
    We first construct a quasi-isomorphism between two 2-term
    complexes of abelian groups:
    \begin{equation*}
    \begin{tikzcd}
        \tfrac{C_1(\W, \hat L)}{B_1(\W, \hat L)}\ar[r, "d"]\ar[d] &
        \hat L\ar[d]\\
        H_1(\W, \hat L)\ar[r, "0"] &
        \hat L_\Gamma.
    \end{tikzcd}
    \end{equation*}
    
    Since $\hat L_0$ is a free abelian group of finite rank, the short
    exact sequence (\ref{hello}) splits, hence there exists a retract
    $\frac{C_1(\W, \hat L)}{B_1(\W, \hat L)}\to H_1(\W, \hat L)$.
    Now take this map as the first vertical map, and let the second
    vertical map be the natural projection $\hat L\to\hat L_\Gamma$.
    It is clear that this defines a quasi-isomorphism.

    Back to the original proposition. It is easy to see that
    the quasi-isomorphism above is continuous \wrt the topology
    on $\frac{C_1(\W, \hat L)}{B_1(\W, \hat L)}$ and $H_1(\W, \hat L)$.
    Hence we have a quasi-isomorphism between two 2-term
    complexes of abelian sheaves:
    \begin{equation*}
    \begin{tikzcd}
        \tfrac{C_1(\W, \hat L)}{B_1(\W, \hat L)}^\circ\ar[r, "d"]\ar[d] &
        \underline{\hat L}\ar[d]\\
        H_1(\W, \hat L)^\circ\ar[r, "0"] &
        \underline{\hat L_\Gamma}.
    \end{tikzcd}
    \end{equation*}
\end{proof}

\begin{prop}
There is a canonical pairing
$\T\times\Locc\to\BGm$.
\label{pairing}
\end{prop}
\begin{proof}
Let $\mathscr P$ be a Picard stack represented by the 2-term complex
$[H\to X]$ where $H$ and $X$ are represented by ind-schemes.
Let $p, a:H\times X\to X$ be the projection and the action map.
A $\Gm$-torsor $\mathcal E$ on $\mathscr P$ is equivalent to a pair
$(\mathcal E_0, \phi)$ where $\mathcal E_0$ is a $\Gm$-torsor on
$X$, and $\phi:p^*\mathcal E_0\to a^*\mathcal E_0$ is an isomorphism
of $\Gm$-torsors such that
\begin{itemize}
    \item the restriction of $\phi$ to $\{e\}\times X$ is the identity.
    \item $\phi$ satisfies the usual cocycle condition.
\end{itemize}

In the special case that $\mathcal E_0$ is the trivial $\Gm$-torsor,
$\phi$ is equivalent to an invertible function $f$ on $H\times X$ such
that
\begin{itemize}
    \item the restriction of $\phi$ to $\{e\}\times X$ is 1.
    \item For $h_1, h_2\in H$ and $x\in X$, we have $f(h_1 h_2, x) =
    f(h_1, h_2 x) f(h_2, x)$.
\end{itemize}

The desired pairing is induced by two pairings:
\begin{alignat*}{3}
\langle\;,\;\rangle: && \hat T & \times \underline{\hat L} && \to \Gm,\\
[\;,\;]: && \ \hZ^1(\W, \hat T)
& \times \tfrac{C_1(\W, \hat L)}{B_1(\W, \hat L)}^\circ
&& \to \Gm.
\end{alignat*}
The first pairing is the natural pairing, since $\hat L$ is the character
lattice of $\hat T$. The second pairing is given by
\begin{alignat*}{2}
	Z^1(\W, \hat T(R)) \,\times\, & C_1(\W, \hat L) && \to R^\times\\
	(\phi,\ &\psi) && \mapsto \sum_{w\in W} \langle
	\phi(w), \psi(w) \rangle.
\end{alignat*}
This pairing factors through $C_1(\W, \hat L)/B_1(\W, \hat L)$ by a
direct computation.

Let $(\phi, x)\in\hZ^1(\W, \hat T)\times\hat L$, and let
$(t, \psi)\in\hat T\times\frac{C_1(\W, \hat L)}{B_1(\W, \hat L)}$.
We define the pairing $\T\times\Locc\to\BGm$ by specifying the
function $f$ as discussed above:
\[
f(\phi, x, t, \psi) = [\phi, \psi] \cdot \langle t, x+d\psi\rangle.
\]
To check the cocycle condition boils down to check that for
any $t, \psi$:
\[
	[dt, \psi] = \langle t,d\psi\rangle
.\]
The verification is straightforward:
\begin{align*}
	[dt, \psi] &= \sum_{w\in W}\langle dt(w), \psi(w)\rangle\\
	&=\sum_{w\in W}\langle wt - t, \psi(w)\rangle\\
	&=\Big\langle t, \sum_{w\in W} w^{-1}\psi(w) - \psi(w)\Big\rangle\\
	&=\langle t, d\psi\rangle. \qedhere
\end{align*}
\end{proof}

\begin{lem}
    Assume that 2 is invertible in $S$. Then
    there is a canonical isomorphism
    \[
            (\rfc{T(F), \Gm})^\vee\cong\BTF.
    \]
\end{lem}
\begin{proof}
    We first claim that for every $i$ there is an isomorphism
    $(\repfun_{T(F)/V_i, \Gm})^\vee\cong\B T(F)/V_i$.
    Recall that the inner Hom between Picard stacks can be
    calculated in the category of abelian sheaves:
    \[
    (\Homm(\mathscr P_1, \mathscr P_2))^\flat =
    \tau_{\leq 0}R\Homm(\mathscr P_1^\flat, \mathscr P_2^\flat).
    \]
    In degree $-1$, we have
    \[
    H^{-1}R\Homm(\repfun_{T(F)/V_i, \Gm}, \BGm) =
    H^{0}R\Homm(\repfun_{T(F)/V_i, \Gm}, \Gm) = T(F)/V_i.
    \]
    In degree $0$, we need to check the vanishing of
    \[
    H^0 R\Homm(\repfun_{T(F)/V_i, \Gm}, \BGm) = 
    H^1 R\Homm(\repfun_{T(F)/V_i, \Gm}, \Gm).
    \]
    Since $T(F)/V_i$ is a finitely generated abelian group,
    $\repfun_{T(F)/V_i, \Gm}$ decompose as a product of
    $\boldsymbol{\mu}_n$ and $\Gm$. It is well known that $\mathscr
    Ext^1(\boldsymbol{\mu}_n, \Gm)=0$ and $\mathscr Ext^1(\Gm, \Gm)=0$,
    hence degree $0$ vanishes. This proves the claim.

    Next, let $\{A_i\}$ be a filtered directed system of complex of
    abelian sheaves, and let $B$ be another such complex. Then there
    is a natural isomorphism
    \[
    R\Homm(\colim A_i, B) \cong \lim R\Homm(A_i, B).
    \]
    Applying the Milnor exact sequence yields a short exact sequence
    \[
    0\to\mathrm{Rlim}^1 H^{n + 1} R\Homm(A_i, B)
    \to H^n R\Homm(\colim A_i, B)
    \to \lim H^n R\Homm(A_i, B)\to 0.
    \]
    In our case, we apply this with $A_i = \repfun_{T(F)/V_i,\Gm}$ 
    and $B=\BGm$.
    For the lemma to hold, it suffice to show the $\mathrm{Rlim}^1$
    terms vanishes for $n = -1$ and $0$. For $n = -1$, the first term
    becomes
    \[
    \mathrm{Rlim}^1 H^0 R\Homm(\repfun_{T(F)/V_i, \Gm}, \BGm)
    = \mathrm{Rlim}^1 \mathscr{E}\!xt^1(\repfun_{T(F)/V_i}, \Gm)=0.
    \]
    For $n = 0$, we need to show the vanishing of
    \[
    \mathrm{Rlim}^1 H^1 R\Homm(\repfun_{T(F)/V_i, \Gm}, \BGm)
    = \mathrm{Rlim}^1 \mathscr{E}\!xt^2(\repfun_{T(F)/V_i}, \Gm).
    \]
    Write $T(F) = M\oplus T(F)_0$, where $M$ is a free abelian
    group of finite rank and $T(F)_0$ is the maximal compact open
    subgroup of $T(F)$. Then
    \[
    \repfun_{T(F)/V_i, \Gm}\cong
    \repfun_{M, \Gm}\oplus\repfun_{T(F)_0/V_i, \Gm},
    \]
    where the first summand is a product of copies of $\Gm$ and the
    second summand is a finite product of $\boldsymbol{\mu}_n$. Since
    the first summand is constant in $i$, its $\mathscr Ext$-groups
    are constant, and contribute nothing to the derived limit. By
    \cite[\S~2, Prop.~1]{Breen}, when 2 is invertible in $S$, we
    have
    \[
    \mathscr{E}xt^2(\boldsymbol{\mu}_n, \Gm) = 0,
    \]
    hence the second summand also contributes nothing to
    $\mathrm{Rlim}^1$.
\end{proof}

\makepoint\textbf{Proof of} $\T\cong\Locc^\vee$\textbf{.}
We form the following diagram where the first row is the split exact
sequence (\ref{tor-ses}), and the second row is the dual of the
split exact sequence (\ref{ses}).
The first vertical map is the isomorphism in the previous lemma,
the second one is induced by the pairing $\T\times\Loc\to\BGm$,
and the third one is the tautological map
$\hat L_\Gamma\to(\hat L_\Gamma^\vee)^\vee$:
\[
\begin{tikzcd}
1\ar[r]& \BTF \ar[r]\ar[d, "\vsim"]&
\T \ar[r]\ar[d]&
\hat L_\Gamma \ar[r]\ar[d, "\vsim"]& 1\\
1\ar[r]& (\rfc{T(F), \Gm})^\vee \ar[r]&
\Locc^\vee \ar[r]&
\Homm(\hat L_\Gamma, \BGm)^\vee \ar[r]& 1.
\end{tikzcd}
\]
Since $\hat L_\Gamma$ is a finitely generated abelian group, it is
dualizable, hence the third vertical map is an isomorphism.

The two visible squares in the diagram commute by construction of the
pairing. Since both rows splits, the square of the +1 maps also
commutes. This concludes the proof that $\T\cong\Locc^\vee$.

\begin{proof}[Proof of theorem \ref{main}.] Now we prove the main
theorem with $\T$ instead of $\Tor$. We established above that for
each torus $T$ the pairing defined in Proposition~\ref{pairing}
induces an isomorphism $\T\cong\Locc^\vee$. It remains to check that
they are a unique family of line bundles on $\T\times\Locc$ satisfying
the three conditions in the theorem.

Condition (a) is equivalent to the following proposition which
one verifies by a direct computation.
\begin{prop}
	Let $g: S\to T$ be a homomorphism between tori. It induces
    $\hat g: \hat T\to\hat S$ and $g^L:\hat L_S\to\hat L_T$.
    Let $f_T(\phi, x, t, \psi)$ and $f_S(\phi, x, t, \psi)$
    be the functions used to construct the Poincar\`e line
    bundles.
    Then for $(\phi, x)\in\hZ^1(\W, \hat T)\times\hat L_S$
    and for $(t, \psi)\in\hat T\times\frac{C_1(\W,\hat L_S)}
    {B_1(\W,\hat L_S)},$, we have
    \begin{equation}
    \label{ftfs}
    f_T(\phi, g^L(x), t, g^L(\psi))
    = f_S(\hat g(\phi), x, \hat g(t), \psi).
    \end{equation}
\end{prop}

Condition (b) follows from the definition of the pairing.

We verify condition (c). Let $S = \Res_{E/F}\Gm$ an induced torus
and $T = \Gm$ a split torus over $E$. The cocharacters of the two torus
satisfy
\[
	\hat L_S = \Ind^{\Gal(E/F)}_{\quad*}\hat L_T
.\]
Notice that by definition,
\[
\Loc_{{}^L S, F}
= \ch\!\left(\tau_{\leq 0}(C^\bullet(\W, \hat S)[1])\right),
\quad
\T_{S/F} = \ch\!\left(\tau_{\geq -1}(C_\bullet(\W, \hat L_S))\right).
\]
We apply Shapiro's lemma to get isomorphisms
\begin{align*}
	\Loc_{{}^L S, F} &= \ch\!\left(
	\tau_{\leq 0}(C^\bullet(\W, \hat S)[1])
	\right)\\
	&\cong \ch\left(
	\tau_{\leq 0}(C^\bullet(\Ex, \hat T)[1])\right) = \Loc_{\cT, E},
\end{align*}
\begin{align*}
	\T_{S/F} &= \ch\!\left(
	\tau_{\geq -1} C_\bullet(\W, \hat L_S)\right)\\
	&\cong \ch\!\left(
	\tau_{\geq -1} C_\bullet(\Ex, \hat L_T)\right)
    = \T_{T/E}.
\end{align*}
The cap product used to define the pairings clearly intertwines
with the Shapiro isomorphisms, and hence condition (c) is satisfied.

By condition (b) and (c), any such family of pairings is uniquely
determined on all induced torus. For every torus $T$, there is always
an induced torus $S$ and a homomorphism $S\to T$ such that the
induced homomorphism on cocharacters $\hat L_S\to\hat L_T$ is
surjective. It suffice to show the pairing for $T$ is determined
by that of $S$. Since $\hat L_S\to\hat L_T$ is surjective,
by (\ref{ftfs}), the function $f_S$ determines the function $f_T$,
hence determines the line bundle for $T$. This concludes the proof of
the main theorem.
\end{proof}

We finish the section with
\begin{prop}
    We have a canonical isomorphism $\T\cong\Tor$.
\end{prop}
\begin{proof}
We start by defining a quasi-isomorphism between two 2-term complexes
of abelian groups
\begin{equation}
\begin{tikzcd}
\frac{C_1(\W, \hat L)}{B_1(\W, \hat L)} \ar[r, "d"]\ar[d, "\cores"]
&\hat L \ar[d, "c'"]\\
T(E) \ar[r, "d"] & \Zalg(\W, T(E)).
\end{tikzcd}
\label{tor-tor-diag}
\end{equation}
The first vertical map is defined to be the corestriction map
\[
	C_1(\W, \hat L)/B_1(\W, \hat L) \xrightarrow{\cores}
	C_1(\Ex, \hat L)/B_1(\Ex, \hat L)\cong T(E).
\]
It sends
$w\otimes x$ to $\sum_{\tau} w_\tau x\otimes\delta(w_\tau, w)$.

Recall that $\Zalg(\W, T(E))$ is defined as a fibre product
\[
        \hat L^\Gamma\times_{Z^1(\Ex, T(E))^\Gamma}Z^1(\W, T(E))
,\]
to define the second vertical map $c'$ amounts to defining two maps
$\hat L \to \hat L^\Gamma$ and $\hat L \to Z^1(\W, T(E))$ with the
usual compatibility condition. We define the first map to be the
norm map, and the second to be the composition
\[
	\hat L \to Z^1(\Ex, T(E)) \xrightarrow{\inf}
	Z^1(\W, T(E))
.\]
Together, they define a map into the fibre product because
$(\res\circ\inf)$ coincides with the norm map.

The inflation map is defined as follows. Let $\phi:\Ex\to
T(E)$ be a cocycle. The inflation of this cocycle is \emph{usually} given by
\[
	g\mapsto \sum_{\tau}w_\tau^{-1}\phi(\delta(w_\tau, g))
.\]
However, we define the inflation map with a \emph{different} system of
representatives given by
$\{\tilde w_\tau = w_{\tau^{-1}}^{-1}, \tau\in\Gamma\}$.
Explicitly, for $x\in\hat L$, $c'(x)$ is the cocycle
(for notation, see section \ref{notation-galois})
\begin{align*}
g&\mapsto \sum_{\tau} \tilde w_\tau^{-1}(x\otimes
\tilde \delta(\tilde w_\tau, g))\\
&= \sum_{\tau} \tilde w_\tau^{-1}x\otimes
\tilde w_\tau^{-1}\tilde \delta(\tilde w_\tau, g)\\
&= \sum_\tau w_{\tau^{-1}} x\otimes \delta(g, w_{g^{-1}\tau^{-1}})
.\end{align*}

Now we verify that the diagram (\ref{tor-tor-diag}) is commutative,
i.e., $c'\circ d = d\circ\cores$.
The cocyle $d\circ\cores(w\otimes x)$ is given by
\begin{equation*}
	g\mapsto \sum_{w_\tau}gw_\tau x\otimes g\delta(w_\tau, w)
	- \sum_{w_\tau}w_\tau x\otimes\delta(w_\tau, w)
.\end{equation*}
We change the variable $\tau\to g^{-1}\tau$ in the first summation,
then the right hand side becomes
\begin{equation}
\label{happy}
\sum_\tau w_\tau x\otimes (g\delta(w_{g^{-1}\tau}, w)
- \delta(w_\tau, w))
=\sum_\tau w_\tau x\otimes (\delta(g, w_{g^{-1}\tau w})
- \delta(g, w_{g^{-1}\tau})).
\end{equation}
On the other hand, $c'\circ d(w\otimes x)$ is the cocycle
\[
	g\mapsto \sum_\tau
	w_{\tau^{-1}}w^{-1}x\otimes\delta(g, w_{g^{-1}\tau^{-1}})
	- \sum_\tau
	w_{\tau^{-1}}x\otimes\delta(g, w_{g^{-1}\tau^{-1}})
.\]
We change the variable $\tau\to w^{-1}\tau$ in the first summation.
The right hand side becomes
\begin{align*}
\sum_\tau w_{\tau^{-1}}x\otimes(\delta(g, w_{g^{-1}\tau^{-1}w})
-\delta(g, w_{g^{-1}\tau^{-1}}))
.\end{align*}
This equals (\ref{happy}).

The homomorphism (\ref{tor-tor-diag}) is a quasi-isomorphism. Indeed,
Corollary \ref{lan-result} shows that $\cores$ induces an
isomorphism on cohomology in degree $-1$, while Corollary
\ref{kot-result} shows that $c'$ induces an isomorphism on cohomology
in degree $0$.

Since the maps defined above are continuous homomorphisms, we have
a quasi-isomorphism between 2-term complexes of abelian sheaves
\begin{equation*}
\begin{tikzcd}
\frac{C_1(\W, \hat L)}{B_1(\W, \hat L)}^\circ
\ar[r, "d"]\ar[d]
&\underline{\hat L} \ar[d]\\
T(E)^\circ \ar[r, "d"] & \Zalg(\W, T(E))^\circ.
\end{tikzcd}
\end{equation*}
Hence, $\T\cong\Tor$.
\end{proof}

\section{Quasi-Coherent Sheaves on \texorpdfstring{$\Tor$}{Tor}}
\label{qcoh-tor}
\noindent
In this section, we show that $\Tor$ plays the role of an integral
version of $\BT$ in the following sense.

\begin{prop}
\label{BT-Tor}
Let $\Lambda = \Zlb, \Qlb$, or $\Flb$, and let $\Shv(-, \Lambda)$
denote the unbounded derived $\infty$-category of $\ell$-adic sheaves
(see Section \ref{l-adic-sheaf-theory} for a precise definition).
Then there is a canonical equivalence
\[
\Shv(\isoc_T, \,\Lambda)\cong\QCoh(\Tor\otimes\Lambda).
\]
\end{prop}

\begin{definition}
Let $H$ be a locally profinite topological group and let $\Lambda$
be a coefficient ring. We define the $\infty$-category of
representations of $H$ with $\Lambda$-coefficients to be the left
completion of $\mathcal D(\Rep(H, \Lambda)^\heartsuit)$ \wrt to
its standard $t$-structure, where $\Rep(H, \Lambda)^\heartsuit$ is the
Grothendieck category of smooth $H$-representations of
$\Lambda$-modules.
\end{definition}

Recall from Section~\ref{abelian-sheaves} that $\uto$ and $\utf$ are
the abelian sheaves associated to the (locally) profinite groups
$T(F)_0$ and $T(F)$:
\[
\uto = \big(T(F)_0\big)^\circ,\quad \utf = T(F)^\circ,
\]
and that $\BTO$ and $\BTF$ are the classifying stacks of $\uto$ and
$\utf$, respectively. For the remainder of this section, we abuse
notation and use the same symbols $\uto$, $\utf$,
$\BTO$ and $\BTF$ to denote their base change to $\Lambda$.

\begin{prop}
\label{qcoh-t0}
There is a canonical equivalence
\[
        \QCoh(\BTO)\cong\Rep(T(F)_0, \Lambda).
\]
\end{prop}
\begin{proof}
    Since $\uto$ is an affine group scheme,
    $\QCoh(\BTO)^\heartsuit$ identifies with the Grothendieck
    category of algebraic $\uto$-representations on
    $\Lambda$-modules, which is equivalent to the abelian category of
    smooth $T(F)_0$-representations on $\Lambda$-modules, denoted
    by $\Rep(T(F)_0, \Lambda)^\heartsuit$.
    
    Since $\uto$ is an affine group scheme, $\BTO$ has affine
    diagonal. It admits an fpqc cover $\Spec\Lambda\to\BTO$. By Lemma
    9.8 (4) in \cite{Zhu25}, $\QCoh(\BTO)$ is the left completion of
    $\mathcal D(\QCoh(\BTO)^\heartsuit)$. On the other hand,
    by definition, $\Rep(T(F)_0, \Lambda)$ is the left completion of $\mathcal D(\Rep(T(F)_0, \Lambda)^\heartsuit)$. The proposition
    follows.
\end{proof}

\begin{prop}
\label{qcoh-t}
There is a canonical equivalence
\[
        \QCoh(\BTF)\cong\Rep(T(F), \Lambda).
\]
\end{prop}
\begin{proof}
    The map $f:\BTO\to\BTF$ is an fpqc-cover (the map $f$ is not
    quasi-compact, but it is an fpqc-cover). By fpqc-descent,
    $\QCoh(\BTF)$ is equivalent to the category of left comodules for
    the comonad $\mathsf T := f_* \circ f^*$ on $\QCoh(\BTO)$. Under
    the equivalence $\QCoh(\BTO)\cong\Rep(T(F)_0, \Lambda)$ in
    Proposition~\ref{qcoh-t0}, this comonad is nothing other than
    $\mathsf T' = \res\circ\ind$ arising from the adjoint pair:
    \[
            \res:\Rep(T(F), \Lambda)\leftrightharpoons
            \Rep(T(F)_0, \Lambda):\ind.
    \]
    Because $\mathsf T'$ exhibit $\Rep(T(F), \Lambda)$ as the
    category of comodules in $\Rep(T(F)_0, \Lambda)$, the proposition
    follows.
\end{proof}

\begin{proof}[Proof of Proposition \ref{BT-Tor}]
    By definition $\isoc_T = LT/\Adsig LT$. From \cite{Kot85}, we know
    \[
    \isoc_T \cong \bigsqcup_{B(T)} [*/T(F)].
    \]
    On the other hand, Remark~\ref{decomposition-tor} gives the decomposition
    \begin{equation*}
    \Tor = \bigsqcup_{B(T)} \BTF,
    \end{equation*}
    with the same index set $B(T)$. Finally, by
    \cite[Proposition~3.51]{Zhu25} and Proposition \ref{qcoh-t}, we have
    \[
        \Shv([*/T(F)], \Lambda)\cong
        \Rep(T(F), \Lambda)\cong
        \QCoh(\BTF\otimes\Lambda).\qedhere
    \]
\end{proof}

\section{Categorical Local Langlands Correspondence for Tori}
\label{transform}
\noindent
Throughout this section, let $\Lambda$ be a commutative ring with
$p$ invertible, and write $\Tor$ and $\Locc$ for their base change
to $\Lambda$.

In section~\ref{duality}, we constructed a canonical family of
Poincar\'e line bundles $\mathcal L_T: \Tor\times\Locc\to \BGm$. The
Fourier--Mukai transforms associated with them are the integral
categorical local Langlands correspondence for tori:

\begin{thm}
    Assume $p$ is invertible in $\Lambda$.
	Let $\pi_1, \pi_2$ be the projections of $\Locc\times\Tor$ to its
	factors. Let $\Nilp$ denote the zero section of the cotangent
    complex of $\Locc$.
	The Fourier--Mukai transform via the line bundle $\mathcal L_T$
    establish an equivalence of $\infty$-categories that preserves
    t-structures:
	\begin{equation*}
		\pi_{1!}(\pi_2^*(-)\otimes\mathcal L_T^{-1}):
		\QCoh(\Tor)\leftrightharpoons
		\IndCoh_{\Nilp}(\Locc):
		\pi_{2*}(\pi_1^*(-)\otimes\mathcal L_T).
	\end{equation*}
\end{thm}

\makepoint[Coherent sheaves on $\Locc$.]
We denote $\rfc{T(F), \Gm}$ by $\uLoc$ and $\rfc{T(F)/V_n, \Gm}$ by
$\uLocn$.  By Corollary \ref{two-ses}, we have splittings
$\Locc = \uLoc \times \bb{B}\hat T^\Gamma$ and $\Locn = \uLocn
\times \B\hat T^\Gamma$.

Notice that $\uLocn = \Spec(\Lambda[T(F)/V_n])$. The inclusion map
$i:\uLocn\to\uLoc^{(n + 1)}$ is given by the ring homomorphism
$\Lambda[T(F)/V_{n + 1}]\to\Lambda[T(F)/V_n]$, which in turn is
induced by the projection $T(F)/V_{n + 1}\to T(F)/V_n$. Since $i$ is a
finite morphism, so is the inclusion $i:\Locn\to\Locc^{(n + 1)}$, and
the $\ast$-pushforward along $i$ sends a coherent sheaf to a coherent
sheaf:
\begin{gather*}
i_*:\Coh(\uLocn)\to\Coh(\uLoc^{(n + 1)}),\\
i_*:\Coh(\Locn)\to\Coh(\Locc^{(n + 1)}).
\end{gather*}
We define $\Coh(\uLoc) = \varinjlim \Coh(\uLocn)$
and let $\IndCoh(\uLoc)$ be its ind-completion.
Similarly, we define $\Coh(\Locc) = \varinjlim \Coh(\Locn)$,
and let $\IndCoh(\Locc)$ be its ind-completion.

\makepoint[Decomposition of categories.]
Under the splittings $\Locc = \uLoc\times \bb{B}\hat T^\Gamma$ and
$\Tor=\hat L_\Gamma\times\BTF$ (by Proposition~\ref{tor-ses}), the
Poincar\'e line bundle $\mathcal L_T$ on $\Locc\times\Tor$ can be
regarded as the structural sheaf over $\uLoc\times\hat L_\Gamma$ with
both a $\hat T^\Gamma$ and a $T(F)$-action, such that $\hat T^\Gamma$
acts by the character in $x\in\hat L_\Gamma$ on the $x$-component, and
$t\in T(F)$ acts on $\uLoc_n\times\hat L_\Gamma$ by $tV_n \in
\Lambda[T(F)/V_n]$.

Since $\Locc = \uLoc\times \bb{B}\hat T^\Gamma$ is a
$\hat T^\Gamma$-gerbe, coherent sheaves on $\Locc$ decompose as
\begin{equation*}
	\Coh(\Locc) = \bigoplus_{\alpha \in X^*(\hat T^\Gamma)}
	\Coh^\alpha(\Locc)
,\end{equation*}
where each subcategory $\Coh^\alpha(\Locc)$ consists of coherent
sheaves on $\Locc$ that has an action of $\hat T^\Gamma$ via the
character $\alpha$. Since the ind-completion functor preserves
colimits, taking ind-completion on both sides gives us
\begin{equation}
	\label{gerbe}
	\IndCoh(\Locc) = \Ind\left(
        \bigoplus_{\alpha \in X^*(\hat T^\Gamma)}
	\Coh^\alpha(\Locc)
    \right)
    \cong \prod_{\alpha \in X^*(\hat T^\Gamma)}
    \IndCoh^\alpha(\Locc)
.\end{equation}

Denote by $\BTF_\beta$ the connected component of $\Tor$ indexed
by $\beta\in B(T) = X^*(\hat T^\Gamma)$. Then we have a decomposition
by connected components
\[
	\QCoh(\Tor) = \prod_{\beta\in X^*(\hat T^\Gamma)}
	\QCoh(\BTF_\beta)
.\]

\makepoint[D\'evissage.] Let us fix $\beta\in B(T)$. 
By abuse of notation, let $\pi_1$ and $\pi_2$ be projections
of $\BTF_\beta\times\Locc$ to its factors:
\[
\begin{tikzcd}
  & \Locc\times\BTF_{\beta} \ar[dl, "\pi_1"'] \ar[dr, "\pi_2"] & \\
  \Locc & & \BTF_{\beta}.
\end{tikzcd}
\]
Let $\mathcal L_\beta$ be the pullback of the Poincar\'e line bundle along
the open embedding $\BTF_\beta\times\Locc \inj \Tor\times\Locc$.
On the one hand, it is easy to check that for $\mathcal F\in\Coh^\alpha(\Locc)$
and $\beta \neq -\alpha$,
\[
	\pi_{2*}(\pi_1^*\,\mathcal F\otimes \mathcal L_\beta) = 0
.\]
On the other hand, for $\mathcal G\in\QCoh(\BTF_{\beta})$,
\(
        \pi_{1!}(\pi_2^*\,\mathcal G\otimes\mathcal{L}_\beta^{-1})
\)
lies in $\IndCoh^{-\beta}(\Locc)$.
Therefore, showing the equivalence in the theorem is reduced to showing the
following equivalence of categories:
\begin{align*}
	\pi_{1!}(\pi_2^*(-)\otimes\mathcal{L}_{-\alpha}^{-1}):
	\QCoh(\BTF_{-\alpha}) \leftrightharpoons
	\IndCoh_{\Nilp}^\alpha(\Locc):
	\pi_{2*}(\pi_1^*(-)\otimes\mathcal L_{-\alpha})
.\end{align*}

We may further reduce our task as follows.
Notice that $\IndCoh_{\Nilp}^\alpha(\Locc)\cong\IndCoh_{\Nilp}(\uLoc)$.
Write $p_1$ and $p_2$ for the projections of
$\uLoc\times \BTF$ to its factors:
\[
\begin{tikzcd}
  & \quad\uLoc \times \BTF \ar[dl, "p_1"'] \ar[dr, "p_2"] & \\
  \uLoc & & \BTF.
\end{tikzcd}
\]
Let $\mathcal L$ be the pullback of $\mathcal L_{-\alpha}$ along the
isomorphism $\uLoc\times\BTF\to\Locc\times\BTF_{-\alpha}$. Showing the
equivalence in the theorem amounts to showing that there is an
equivalence of categories:
\[
    p_{1!}(p_2^*(-)\otimes\mathcal L^{-1}):
	\QCoh(\BTF) \leftrightharpoons
	\IndCoh_{\Nilp}(\uLoc):
	p_{2*}(p_1^*(-)\otimes\mathcal L).
\]

\makepoint\textbf{The equivalence.}
In this section we show that the functor
\begin{equation}
p_{2*}(p_1^*(-)\otimes\mathcal L):
\IndCoh_{\Nilp}(\uLoc)\to\QCoh(\BTF)
\label{indcoh-to-qcoh}
.\end{equation}
is an equivalence. Since taking colimit commutes with ind-completion,
we have
\[
\IndCoh_{\Nilp}(\uLoc) = \varinjlim\IndCoh_{\Nilp}(\uLocn)
= \varinjlim\QCoh(\uLocn)
.\]
Let $N\in\QCoh(\uLocn)^\heartsuit$, or equivalently,
let $N$ be a $\Lambda[T(F)/V_n]$-module.
The quasi-coherent sheaf $p_1^*\,N\otimes\mathcal L$
is supported on $\uLocn\times\BTF$. It can be regarded as
the same $\Lambda[T(F)/V_n]$-module $N$ with a compatible $T(F)$-action,
such that the action of $t\in T(F)$ is identified with the
multiplication by $tV_n\in\Lambda[T(F)/V_n]$. The pushforward by $p_{2*}$
simply forgets the $\Lambda[T(F)/V_n]$-module structure.

In summary, the pull-push factors as
\[
        \begin{tikzcd}[column sep=small, cramped]
        \QCoh(\uLocn)^\heartsuit
        \ar[rr, "\pi_{2*}(\pi_1^*(-)\otimes\mathcal L)"]
        \ar[rd, "\sim"]&&
        \QCoh(\BTF)^\heartsuit.\\
        &\Lambda[T(F)/V_n]\mathrm{-Mod} \ar[ru]&
        \end{tikzcd}
\] 
\begin{lemma}
Assume $1/p\in\Lambda$, we have
\[
\varinjlim\Rep(T(F)/V_n, \Lambda) \cong \Rep(T(F), \Lambda).
\]
\end{lemma}
\begin{proof}
We claim that both $\infty$-categories $\Rep(T(F)/V_n, \Lambda)$
and $\Rep(T(F), \Lambda)$ are the derived category of their hearts.

Since $\Rep(T(F)/V_n, \Lambda)$ is equivalent to the category of
$\Lambda[T(F)/V_n]$-modules, it is the derived category of its heart.

It is well known that if a Grothendieck abelian category
$\mathcal C^\heartsuit$ has a set of generators $\{c_i\}_i$
such that $\Ext^\bullet_{\mathcal C^\heartsuit}(c_i, -)$
has finite cohomological dimension, then its derived category
$\mathcal D(\mathcal C^\heartsuit)$ is left complete.

Consider the set of generators $\{\delta_K\}_K$ for all
pro-$p$ open compact subgroup $K$ of $T(F)$. Each $\delta_K =
c\text{-}\!\ind^{T(F)}_K \mathbf{1}$ is a projective object in
$\Rep(T(F), \Lambda)^\heartsuit$, where $\mathbf{1}$ denotes the
trivial representation of $K$. Therefore $\Rep(T(F), \Lambda)$
is also the derived category of its heart.

When $V_n$ is pro-$p$ (this is true for sufficiently large $n$), the
inclusion
\[
\Rep(T(F)/V_n, \Lambda)\to\Rep(T(F), \Lambda)
\]
is fully faithful. This is because the inclusion
$\Rep(T(F)/V_n, \Lambda)^\heartsuit \to\Rep(T(F), \Lambda)^\heartsuit$
is exact, fully-faithful, and sends projective objects to projective
objects ($\delta_{V_n} = \Lambda[T(F)/V_n]$). The same reasoning
shows that
\[
\Rep(T(F)/V_n, \Lambda)\to\Rep(T(F)/V_{n + 1}, \Lambda)
\]
is also fully faithful. Therefore we have a fully-faithful functor
\[
\varinjlim\Rep(T(F)/V_n, \Lambda)\to\Rep(T(F), \Lambda).
\]
This is an equivalence because for any smooth $T(F)$-representation
$M$, $M = \varinjlim M^{V_n}$.
\end{proof}

Therefore, we have
$\IndCoh_{\Nilp}(\uLoc) = \varinjlim\QCoh(\uLocn)
\cong \varinjlim\Rep(T(F)/V_n, \Lambda)
\cong \Rep(T(F), \Lambda) \cong \QCoh(\BTF)$.

\makepoint\textbf{The inverse functor.}
The functor $p_1^*$ preserves limits and is therefore a right adjoint.
We denote its left adjoint by $p_{1!}$. In fact, $p_{1!}$ is the
derived functor of taking $T(F)$-coinvariants. We want to show the
following functor is an equivalence:
\[
p_{1!}(p_2^*(-)\otimes\mathcal L^{-1}):
\QCoh(\BTF)\to\IndCoh_{\Nilp}(\uLoc).
\]

Let $M\in\QCoh(\bb{B}(T(F)/V_n))$. The quasi-coherent sheaf
\(
p_2^*\,M\otimes\mathcal L^{-1}
= M\otimes\Lambda[T(F)/V_n]
\)
is supported on $\uLocn\times\BTF$. Let
$m\otimes f\in M\otimes\Lambda[T(F)/V_n]$, then the action of $t\in T(F)$
is given by
\[
t\cdot(m\otimes f) = (t\cdot m) \otimes (t^{-1}V_n\cdot f).
\]
Due to the fact that $\Lambda[T(F)/V_n]$ is a projective $T(F)/V_n$-module,
we have
\[
	p_{1!}(p_2^*\,M\otimes\mathcal L^{-1})|_{\uLocn} =
	p_{1!}(M\otimes \Lambda[T(F)/V_n])\cong M
,\] 
where the last term $M$ is not the original $T(F)$-module but an
$\Lambda[T(F)/V_n]$-module. In general, take $M = \varinjlim M_n$ where
$M_n = M^{V_n}$. As $p_{1!}$ commute with filtered colimits, we have
\[
        p_{1!}(p_2^*\,M\otimes\mathcal L^{-1}) = M.
\]
In other words, the functor
$\pi_{1!}(\pi_2^*(-)\otimes\mathcal L^{-1})$ is the inverse of
$\pi_{2*}(\pi_1^*(-)\otimes\mathcal L)$.

\makepoint\textbf{The Whittaker condition.} In this subsection, we
show that, under our Langlands correspondence, the regular function
ring corresponds to the Whittaker sheaf.

Let $K$ be a pro-$p$ open compact subgroup of $T(F)$ and let
$\mu$ be the (left) Haar measure on $T(F)$ such that $\mu(K) = 1$.
We define the Hecke algebra of $T(F)$ to be $\Lambda$-valued compactly
supported smooth function on $T(F)$ equipped with the convolution
product \wrt $\mu$ and denote it by $\mathcal H(T(F), \Lambda)$.

The Whittaker sheaf is the Hecke algebra with $T(F)$ acts
by translation.

\begin{prop}
Define the regular function ring $\mathcal O_{\uLoc}$ to consist
of all functions that is supported on $\uLocn$ for some $n$. Then
there is an isomorphism $\mathcal O_{\uLoc}\cong
\mathcal H(T(F), \Lambda)$.
\end{prop}
\begin{proof}
	$G_n = T(F) / V_n$ is a tower of finitely generated abelian
	groups with projections $p_n: G_{n + 1} \to G_{n}$. We have
    $\uLocn = \Spec(\Lambda[G_n])$. There exists
    an integer $N$ such that $V_N\subset K$.

    For $n \geq N$, consider the ring map $\Lambda[G_n]\to\mathcal H(G_n,
    \Lambda)$ that sends an element $g\in G_n$ to
    $\frac{\delta_g}{\mu(V_n)}$. By direct computation, this map is an
    algebra isomorphism (with $\bb{C}$-coefficient, this is the familiar
    discrete Fourier transform).

	The Hecke algebra $\mathcal H(T(F), \Lambda)$ is the union of all
    $\mathcal H(G_n, \Lambda)$, where the inclusion map $\mathcal
    H(G_n, \Lambda)\to \mathcal H(G_{n + 1}, \Lambda)$ is the pullback of
    functions along projections $p_n$. For $n \geq N$, let
    $i_n:\Lambda[G_n]\to\Lambda[G_{n + 1}]$ be the unique map such that
    the following diagram commutes:
    \[
	\begin{tikzcd}
	\Lambda[G_n] \ar[r, "\sim"] \ar[d, "i_n"] &
	\mathcal H(G_n, \Lambda) \ar[d, hook]\\
	\Lambda[G_{n + 1}] \ar[r, "\sim"] &
	\mathcal H(G_{n + 1}, \Lambda).
	\end{tikzcd}
	\]
    It remains to identify the map $i_n$ with the inclusion
	$\mathcal O_{\uLocn}\subset\mathcal O_{\uLoc^{(n + 1)}}$.
    This is clear, since the image of $1$ under $i_n$ is
    an idempotent
    \[
    i_n(1) = \frac{|G_n|}{|G_{n + 1}|}
    \sum_{g\in\operatorname{ker}(p_n)} g,
    \]
    and that $i_n(1)\cdot \Lambda[G_{n+1}]$ is precisely the image of $\Lambda[G_n]$.
\end{proof}

\section{Notation II}
\label{notation-ii}
\refmakepoint[The inertia group.]
\label{inertia-group}
Let $I \subset \Gamma$ be the inertia group of the extension $\F/F$.
Choose an algebraic closure $\overline{F}$ of $F$, and let $\breve
E\subset \overline{F}$ and $\breve F\subset \overline{F}$ denote the
maximal unramified extensions of $E$ and $F$. Write $I_E = 
\Gal(\overline{F}/\Omega_E)$ and $I_F = \Gal(\overline{F}/\Omega_F)$
for the absolute inertia groups of $E$ and $F$.

We define the relative inertia group as $\I := I_F/[I_E, I_E]$.
It is a group extension that corresponds to the canonical class
$u_{\breve E/\breve F}\in H^2(I, I_E^\ab)$ defined in\cite[\S~2.5]{Ser61}:
\[
	1\to I_E^\ab\to\I\to I\to 1
.\]
The triple $(I, I_E^\ab, u_{\breve E/\breve F})$ satisfies the
Tate-Nakayama criterion (see Proposition \ref{triples}).

\makepoint \textbf{Strongly continuous representations.} Let $G$ be a
locally profinite group and let $H$ be a flat affine group scheme of finite type
over $\Zl$. For any $\Zl$-algebra $R$, we say a representation
$f: G\to H(R)$ is \emph{strongly continuous} if for one (and therefore every)
faithful representation $H\to\mathrm{GL}(M)$ on a finite free $\Zl$-module
$M$, for every $m\in M\otimes R$ and every open compact subgroup $G_0$
of $G$, the $\Zl$-module $N\subset M\otimes R$ generated by $f(G_0) m$
is finitely generated, and the resulting representation of $f(G_0)$ on
$N$ is continuous (with $N$ equipped with the $\ell$-adic topology).
We denote the set of strongly continuous representations $f: G\to
H(R)$ by $\scHom(G, H(R))$.

We define the representation functor of strongly continuous representations,
denoted by $\repfunsc_{G, H}$, to be the functor
\[
\repfunsc_{G, H}: \textbf{CAlg}_{\Zl}\to\textbf{Set},\quad
R\mapsto \scHom(G, H(R)).
\] 

\makepoint \textbf{The loop group.} We introduce the notations for
$p$-adic loop groups following \cite{Zhu14}. We denote by $\kappa_F$
and $\mathcal O_F$ the residue field and ring of integers of $F$. Let
$\varpi$ denote a uniformizer of $\mathcal O_F$. For a
$\kappa_F$-algebra $R$, let $W(R)$ denote its ring of Witt vectors,
and let $W_{\mathcal O_F}(R) = W(R)\otimes_{W(\kappa_F)}\mathcal O_F$,
$W_{\mathcal O_F, n}(R) = W(R)\otimes_{W(\kappa_F)}\mathcal
O_F/\varpi^n$.

Let $\mathcal Y$ be a finite-type $\mathcal O_F$-scheme. According to
Greenberg, the following two presheaves on affine $\kappa_F$-schemes:
\[
	L^+_p\mathcal Y(R) = \mathcal Y(W_{\mathcal O_F}(R)),\quad
	L^n_p\mathcal Y(R) = \mathcal Y(W_{\mathcal O_F, n}(R)),
\] 
are represented by schemes over $\kappa_F$. We denote their perfection by
\[
	L^+\mathcal Y = (L^+_p\mathcal Y)^{p^{-\infty}},\quad
	L^n\mathcal Y = (L^n_p\mathcal Y)^{p^{-\infty}},
\] 
and call them $p$-adic jet spaces.

Let $Y$ be an affine scheme of finite type over $F$. The $p$-adic loop
space $LY$ of $Y$ is a perfect space defined by assigning a perfect
$\kappa_F$-algebra $R$ the set
\[
	LY(R) = Y(W_{\mathcal O_F}(R)[1/p])
.\] 
$LY$ is represented by an ind-perfect scheme.

Assume $\mathcal Y$ is an affine scheme of finite type
over $\mathcal O_F$ and $Y = \mathcal Y\otimes_{\mathcal O_F} F$, then
$L^+\mathcal Y\subset LY$ is a closed subscheme.

In \parti, let $\mathfrak T$ be the connected N\'eron model of $T$. We
denote $L^+\mathfrak T$ simply by $L^+T$. Define $\LnT =
\ker(L^+\mathfrak T\to L^n\mathfrak T)$ to be congruence subgroups of
$L^+T$.

\refmakepoint[$\ell$-adic sheaf theory.] \label{l-adic-sheaf-theory}
Let $\Lambda=\Zlb, \Qlb$, or $\Flb$. For
$X\in\AlgStk^\perf_{\kappa_F}$, the stable
$\infty$-category of $\ell$-adic sheaves $\Shv(X, \Lambda)$ is
defined as follows. When the coefficient ring is clear from context,
we simply write $\Shv(X)$.

For an affine scheme $X$ perfectly of finite type over $\kappa_F$, we
denote by $\mathcal{D}^b_c(X, \Lambda)$ the bounded derived
$\infty$-category of $\ell$-adic sheaves on $X$, and by $\Shv(X,
\Lambda)$ its ind-completion.

For an arbitrary perfect affine scheme $X$, we write
$X$ as a filtered limit $X = \lim_\alpha X_\alpha$ with each $X_\alpha$
affine and perfectly of finite type. We denote by $\Shv(X, \Lambda)$
the colimit $\colim_\alpha\Shv(X_\alpha, \Lambda)$ taken \wrt
$!$-pullbacks. It is easy to see the resulting category is independent
of the presentation.

For an arbitrary perfect algebraic stack $\mathcal X$, we denote by
$\Shv(X, \Lambda)$ the limit $\lim_{X\to\mathcal X}\Shv(X, \Lambda)$,
taken over all morphisms $X\to\mathcal X$, where $X$ is a perfect
affine scheme.

See \cite[\S~10]{Zhu25} for a comprehensive development of this
$\ell$-adic sheaf theory and its six-functor formalism.

\makepoint \textbf{Character sheaves and Serre's fundamental group.}
Let $H$ be a connected pro-algebraic group over $\overline{\bb{F}}_p$,
and let $m$ be the multiplication map of $H$. Recall that a character
sheaf on $H$ (with coefficient in $\Lambda$) is a rank-one
$\Lambda$-local system $\Ch_\xi$ on $H$ equipped with an isomorphism
\[
	m^*\,\Ch_\xi\cong\Ch_\xi\boxtimes_\Lambda\Ch_\xi,
\]
which satisfies the usual cocycle conditions.

Denote by $\mathcal P$ the category of commutative pro-algebraic
groups over $k = \overline{\mathbb{F}}_p$ and by $\mathcal P_0$
the category of abelian profinite groups.
Serre defined the fundamental groups $\pi_1$ as the first derived
functor of the right-exact functor $\pi_0:\mathcal P\to\mathcal P_0$,
$\pi_0(G):=G/G^\circ$ (\cite{Ser61}, see also \cite[\S~4.1.1]{Zhu25}
for a different approach).

Assume $H = H_1\times H_2$ is abelian and connected. Then we have
\[
        \pi_1(H) = \pi_1(H_1)\times\pi_1(H_2)
.\]

Assume $H$ is abelian and connected. A key property of this
fundamental group is that the abelian group of character sheaves on
$H$, which we denote by $\mathrm{CS}(H, \Lambda)$, is isomorphic to
the abelian group of continuous rank-1 representations of $\pi_1(H)$:
\[
	\cHom(\pi_1(H), \lx)\cong
	\mathrm{CS}(H, \Lambda)
.\]

In section \ref{geom-llc}, we categorify this
isomorphism into a fully faithful functor (take $H = L^+T$):
\[
	\Ch:\IndCoh\left(\rfsc{\piLT, \Gm}\right)
	\to
	\Shv(L^+T, \Lambda)
.\]
The following theorem, which can be called the geometric class field
theory, will play a key role in the construction of this functor.

\begin{thm}[\cite{Ser61, DW23}]
    Let $T$ be a torus over $F$ and splits over a finite Galois
    extension $E/F$. Let $I$ be the inertia group of the field
    extension $F/E$ and let $I_E^\ab$ be the abelianization of the
    (absolute) inertia group of $E$. Then there is a canonical
    isomorphism
    \[
	\pi_1(L^+T) \cong (\hat L\otimes I^\ab_E)^I
    .\] 
\end{thm}

\makepoint[The fixed-point stack.]
Let $Y$ be an arbitrary stack equipped with an automorphism $\sigma_Y:
Y\to Y$. We define the fixed-point stack of $\sigma_Y$ to be
\[
	\Lsig Y = Y_{\,\,\Gamma,\, Y\times Y,\, \Delta}\, Y
,\] 
where $\Gamma = (1, \sigma_Y)$ is the graph of $\sigma_Y$ and $\Delta$ is
the diagonal embedding. The natural projection to the first factor
defines a canonical morphism
\[
\Lsig Y\to Y
.\]

\begin{prop}
There is a canonical isomorphism of stacks
$\mathcal L_\sigma \BLT \cong \dfrac{LT}{\Adsig LT}$.
\end{prop}
\begin{proof}
By definition $\mathcal L_\sigma \BLT$ is the fibre product
\[
\begin{tikzcd}
\mathcal L_\sigma \BLT \ar[r] \ar[d] & \BLT \ar[d, "\Delta"] \\
\BLT \ar[r, "\mathrm{id}\times\sigma"] & \BLT\times \BLT.
\end{tikzcd}
\]
In other words, for a perfect $\kappa_F$-algebra $R$, $\mathcal
L_\sigma \BLT(R)$ classifies $(\mathcal E, \phi)$ where $\mathcal E$
is an $LT$-torsor over $\Spec R$ and $\phi:\mathcal E \to
\sigma_R^*\mathcal E$ is an $LT$-torsor isomorphism. From this
description, it is straightforward to see that $\mathcal L_\sigma \BLT
\cong \frac{LT}{\Ad_\sigma LT}$.
\end{proof}

\section{Stack of Inertial Langlands Parameters}
\label{rep-stack-of-inertia}

\noindent 
From now on, let $S = \Spec(\Zl)$ be the base scheme. 

\begin{definition}
Let $p: I_F\to\Gamma$ be the natural projection $I_F\to I$ composed
with the inclusion $I\subset\Gamma$. We define the framed functor of
inertial Langlands parameters to be
\[
\Y: \textbf{CAlg}_{\,\Zl} \to \textbf{Set}, \quad 
R\mapsto \repfunsc_{I_F, \cT}\times_{\repfunsc_{I_F, \Gamma}} \{p\}.
\]
It carries a natural action of $\hat T$, and the stack of
inertial Langlands parameters is the quotient stack
$\X = [\Y/\hat T]$.
\end{definition}

We will see in Section~\ref{X-ind-algebraic} that $\X$ is an
ind-algebraic stack.

\begin{prop}
We have a canonically split short exact sequence
\begin{equation}
	\label{ses3}
	1 \to\Homm(\hat L_I, \BGm)\to\X\to\rfsc{\piLTvalue, \Gm}\to 1
.\end{equation}
\end{prop}
\begin{proof}
Since $I_E^\ab$ is profinite, $(I, I_E^\ab, u_{\breve E/\breve F})$
satisfies the Tate-Nakayama criterion, and $H_1(\I, \hat L) = (\hat
L\otimes I^\ab_E)^I$, all the results of Section~\ref{lemmas-in-HA}
and \ref{short-exact-sequence} applies. To prove the proposition, it
suffices to prove a slightly modified version of Lemma~\ref{cont-iff}
as follows.
\end{proof}
\begin{lemma}
Let $A$ be a locally profinite abelian group, let $\Gamma$ be a
finite group acting on $A$ and let $G$ be a group extension of $A$ and
$\Gamma$. Let $L$ be a free abelian group of finite rank with a
$\Gamma$-action, and let $\hat L$ be its dual.

For any $\Zl$-algebra $R$, we define $\pi$ to be
\[
\begin{tikzcd}
	\bar Z^1(G, L\otimes\Rx) \ar[r, "\pi"] &
	\Hom(H_1(G, \hat L), R^\times)\\
	\phi:G \to L\otimes R^\times \ar[r, mapsto] &
	\Big(
	\sum g\otimes x \mapsto
	\sum \langle \phi(g), x \rangle
	\Big).
\end{tikzcd}
\]
Then
$\phi\in\bar Z^1(G, L\otimes R^\times)$ is strongly continuous if and
only if $\pi(\phi)$ is strongly continuous.
\end{lemma}
\begin{proof}
Suppose $R$ is a $\Zl$-algebra. We want to show
$\psi:H_1(G, \hat L)\to R^\times$ is strongly continuous if and only if
$\psi\circ\cores$ is.

By definition, a representation $f:H\to R^\times$ is
strongly continuous if the $\Zl$-module $M$ generated by $f(H)$ is
finite over $\Zl$, and the action of $H$ on $M$ is continuous with $M$
equipped with the $\ell$-adic topology. Equivalently, it is to require
that the map $f:H\to M$ is continuous. Since $\cores$ has finite
cokernel, $\Im(\psi)/\Im(\psi\circ\cores)$ is finite, hence
$\Im(\psi)$ generates a finite $\Zl$-module if and only if
$\Im(\psi\circ\cores)$ does. Let $M$ be the $\Zl$-module generated by
$\Im(\psi)$. Clearly $\psi:H_1(G, \hat L)\to M$ is continuous if and
only if $\psi\circ\cores:H_1(A, \hat L)\to M$ is continuous.

The remainder of the proof follows exactly as in Lemma \ref{cont-iff}.
\end{proof}

\refmakepoint $\X$ \textbf{as an ind-algebraic stack.}
\label{X-ind-algebraic}
To understand the geometry of $\rfsc{\piLTvalue, \Gm}$,
we first observe that the abelian group $I^\ab_E$ canonically
splits into its tame and wild parts as a $\Gamma$-module. By class
field theory, we have
\[
    I^\ab_{\F} = \varprojlim_{\substack{
    K/\F\, \text{finite Galois},\\
    \text{unramified}
}}
\mathcal O_K^\times.
\]
Each term $\mathcal O_K^\times = \mu_K\times U^{(1)}_K$ splits as a
$\Gamma$-module, where $\mu_K$ is the group of roots of unity in $K$
and $U^{(1)}_K = 1 + \mathfrak{m}_K$ is the group of principal units.
This splitting clearly respects the norm map, which is the transition
map in the projective limit. Denote by $I^\tame_E$ (resp. $P$) the
inverse limit of the groups $\mu_K$ (resp. $U^{(1)}_K$). Then we have
a splitting
\[
I^\ab_E = I^\tame_E \times P
.\]

Using the splitting of the group $I^\ab_E$, we obtain
\[
\rfsc{\piLTvalue, \Gm} =
\rfsc{(\hat L\otimes I^\tame_E)^I, \Gm}\times
\rfsc{(\hat L\otimes P)^I, \Gm}
\]
which is characterized by that $(\hat L\otimes I^\tame_E)^I$ is
a prime-to-$p$ profinite group and $(\hat L\otimes P)^I$
is a pro-$p$ group. This gives us
\begin{equation}
\label{decomp_X}
\rfsc{\piLTvalue, \Gm} =
\bigsqcup_{x\in\Xi} \;\rfsc{\hat J, \Gm}.
\end{equation}
where
$\Xi = \rfsc{(\hat L\otimes P)^I, \Gm}$, $x\in\Xi$ is a closed
point and $\hat J = (\hat L\otimes I^\tame_E)^I$.
\smallskip

The functor $\rfsc{\hat J, \Gm}$ can be made explicit by
embedding it into an algebraic torus. We choose $\tau$ a topological
generator of $I^\tame_E$, and write $J$ for the finitely generated
abelian group $(\hat L\otimes \langle\tau\rangle)^I\subset\hat J$.
Denote by $\hat H$ the torus $\repfun_{J, \Gm}$.
There is an inclusion $\rfsc{\hat J, \Gm}\to\hat H$ by
restricting a representation of $\hat J$ to $J$.

Let $\hat H(\Flb)^p$ denote the subset of $\hat H(\Flb)$ consisting
of all elements of order prime to $p$. Denote by $\hat H^p
\subset\hat H$ the subfunctor which is the union of all closed
subschemes 
\[
i_Z:Z\subset\hat H
\]
that are finite over $\Zl$ and satisfy $Z(\Flb)\subset\hat H(\Flb)^p$.
Note that by construction $\hat H^p$ has a natural structure of an
ind-scheme.

\begin{prop}
The inclusion $\rfsc{\hat J, \Gm}\to\hat H$ identifies
$\rfsc{\hat J, \Gm}$ with $\hat H^p$.  Consequently,
$\rfsc{\hat J, \Gm}$ is an ind-scheme.
\label{embedding}
\end{prop}
\begin{proof}
Let $A$ be a $\Zl$-algebra.
Every $\phi\in\mathcal R_{\hat J, \Gm}(A)$
is a strongly continuous homomorphism
$\phi:\hat J\to A^\times\subset A$, or equivalently,
the image of $\varphi$ in $A$ spans a finitely generated
$\Zl$-module $M$, and the map $\phi:\hat J\to M$ is continuous.

The image of $\varphi\in\hat H$ gives a ring homomorphism
which we also denote by $\varphi:\Zl[J]\to A$.
Let $\varphi$ factor as $\Zl[J]\to R\inj A$,
where $\Spec R$ is the schematic image of $\Spec A$.
Our conditions on $\varphi$ imply that $R$ is a finite $\Zl$-algebra.

Consequently, every $\Fl$-points of $\Spec R$ lifts to an
$\Fl$-point of $\Spec A$, and its image in $\hat H(\Fl)$
has finite and prime-to-$p$ order, because $\varphi:
\Zl[J]\to A$ can be lifted to a continuous map
$\Zl[\hat J]\to A$.
\end{proof}

\begin{cor}
    $\X$ is an ind-algebraic stack.
\end{cor}

\begin{cor}
    \label{component-in-hat-h}
    Let $\Lambda = \Flb$ or $\Qlb$, and let $r = \rank\hat T$. Then
    each closed point $x\in\rfsc{\hat J, \Gm}\subset\hat H$ lies on a
    component $X$ that is isomorphic to $\Spf\Lambda[[X_1, \cdots,
    X_r]]\times\B\hat T^I$. Consequently, $\Sing(\X) =
    \X\times\B(\mathfrak t^*)_I$, where $\mathfrak t^*$ denotes the
    dual of the Lie algebra of $\hat T$. Let $0$ denote the zero
    section of $\Sing(\X)$. We have $\Coh_0(\X) = \Coh(\X)$ and
    $\IndCoh_0(\X) = \IndCoh(\X)$.
\end{cor}

\begin{cor}
    \label{generator-for-coh} Let $\Lambda = \Flb$ or $\Qlb$. The
    category $\Coh(\X\otimes\Lambda)$, and hence
    $\IndCoh(\X\otimes\Lambda)$, is generated by skyscraper sheaves
    $\mcal O_x$ for all closed points $x\in\X\otimes\Lambda$.
\end{cor}

Let us fix a Frobenius element $\sigma$. The Frobenius action on $\X$
sends a cocycle $c: I_F\to \hat T$ to $\sigma c$ such that $\sigma
c(g) = {}^\sigma c(\sigma^{-1} g \sigma)$.

\begin{prop}
\label{fixed-point-stack}
The fixed-point stack $\Lsig\X$ is canonically isomorphic to
$\Locc\otimes\Z_\ell$.
\end{prop}
\begin{proof}
There is a natural map $\Locc\otimes\Z_\ell\to\Lsig\X$ induced by
restricting a cocycle $c:W_F\to\hat T$ to $c':I_F\to\hat T$. The goal
is to show that this map is surjective on objects and is a bijection
on the automorphism group of an object.

Objects of $\Lsig\X$ can be represented by a tuple $([c], t)$ for
$c\in Z^1(I_F, \hat T)$ and $t\in\hat T$ such that $\sigma c =
c + dt$.  It is the image of $\tilde c\in Z^1(W_F, \hat T)$ if
$\tilde c$ satisfies $\tilde c|_{I_F} = c$ and $\tilde c(\sigma)=t$.

To show that such $\tilde c$ exists, we first need to verify that $c$
factors through a finite quotient of $I_F$. The restriction of $c$ to
the wild inertia $P_F$ clearly factors through a finite quotient.
Therefore, it suffices to show the claim for the restriction of $c$ to
$I_F^\tame$.

Let the action of $I_F^\tame$ on $\hat T$ factor through an open subgroup $K$.
Then for all $g\in K$, we have
$c(g) = {}^\sigma c(\sigma^{-1}g\sigma)$. Henceforth ${}^\sigma c(g) =
c(\sigma g\sigma^{-1}) = c(g^q) = c(g)^q$. Let $N$ be a positive
integer so that $\sigma^N$ acts trivially on $\hat T$, we have $c(g) =
c(g)^{q^N}$. That is, for every $g\in K$, $c(g^{q^N - 1}) =
1$. We conclude that $c$ factors through an open subgroup of $I_F^\tame$ as
$\{g^{q^{N - 1}} | g\in K\}$ is open. 

We define $\tilde c(\sigma^n)$ 
recursively via
\[
	\tilde c(\sigma^n) = t + {}^\sigma\tilde c(\sigma^{n-1})
\] 
for all $n\in\mathbb N$ and define
\begin{align*}
	&\tilde c(\sigma^{-n}) = -
	{}^{\sigma^{-n}}\tilde c(\sigma^n),\\
	&\tilde c(\sigma^m g) = \tilde c(\sigma^m) +
	{}^{\sigma^{n}} c(g),
\end{align*}
for all $n\in\mathbb N, m\in\mathbb Z$ and for all $g\in I_F$. The
condition $\sigma c = c + dt$ ensures that $\tilde c$ is a 1-cocycle.

An automorphism of the object $([c], t)$ is represented by an
element $w\in\hat T$ such that $dw = 0$ and the following diagram
commutes:
\[
\begin{tikzcd}
 c \ar[r, "t"]\ar[d, "w"] & \sigma  c \ar[d, "\sigma w"]\\
 c \ar[r, "t"] & \sigma c.
\end{tikzcd}
\] That is, $w\in\hat T^{I_F}$ and $\sigma w = w$, hence
$w\in\hat T^{W_F}$. This automorphism group is exactly in bijection
with the automorphism group of an object of $\Locc$.
\end{proof}

\begin{cor}
    \label{image-of-i}
    Let $\Lambda = \Flb$ or $\Qlb$. Let
    \[
    i:\Locc\otimes\Lambda\cong\Lsig\X\otimes\Lambda
    \to\X\otimes\Lambda
    \]
    be the isomorphism in the proposition composed with the
    canonical morphism of the fixed-point stack.
    Then the image of $\Locc(\Lambda)\to\X(\Lambda)$ is canonically
    isomorphic to the set of $\sigma$-invariant characters
    $\cHom(\piLT, \lx)^\sigma$.
\end{cor}
\begin{proof}
    From Proposition~\ref{embedding}, we know that the set of
    $\Lambda$-points of $\X$ can be identified with a subset
    of $\cHom(\piLT, \lx)$. Moreover, this identification respects
    the $\sigma$-action. It is
    clear that the image of $\Locc(\Lambda)\to\X(\Lambda)$ is
    contained in $\cHom(\piLT, \lx)^\sigma$. For every
    character $\chi\in\cHom(\piLT, \lx)^\sigma$, the first
    half of the proof of Proposition~\ref{fixed-point-stack}
    precisely shows it lifts to some $\Lambda$-point of $\Locc$.
\end{proof}

\begin{lemma}
    \label{invariant-char}
    There is an isomorphism
    \begin{equation*}
    \cHom(\piLT, \lx)^\sigma\cong\cHom(T(F)_0, \lx).
    \end{equation*}
\end{lemma}
\begin{proof}
    Using Lemma 2.7 in \cite{DW23}, we have
    \begin{align*}
            &\cHom(\piLT, \lx)^\sigma\\
            =& \cHom(\piLT_\sigma, \lx)\\
            =& \cHom(L^+T(\kappa_F), \lx)\\
            =& \cHom(T(F)_0, \lx).\qedhere
    \end{align*}
\end{proof}

\section{Categorical Trace and Convolution Pattern}
\label{cat-trace}\noindent
The reader may skip this section on a first reading. One should,
however, take note of Section \ref{point-def-psi} and Corollary
\ref{cor-def-phi} before moving on, as they form part of the basic
input for what follows.

In this section we first recall the definition of the categorical
trace and introduce a convolution pattern that frequently arises in
the geometric Langlands program. The main statement is Proposition
\ref{cor860}, which describes the computation of categorical
trace and show that the trace of a sheaf category can often be
identified with a full subcategory of sheaves on a suitable
fixed-point object.

\makepoint[Hochschild homology.]
Let $\mathcal R$ be a symmetric monoidal category, $A$ and $B$ be two
associative algebras in $\mathcal R$. We denote by $A^\rev$ (resp.
$B^\rev$) the algebra $A$ (resp. $B$) with the reversed
multiplication, then an $A$-$B$-bimodule can be regarded as a left
$(A\otimes B^\rev)$-module or a right $(B\otimes A^\rev)$-module. For
an $A$-$A$-bimodule $F$, the Hochschild homology of $F$, if it exists,
is defined as
\[
\Tr(A, F) = A\otimes_{A\otimes A^\rev} F\in\mathcal R.
\]
However, the Hochschild complex of $F$ always exists. It is
given by
\[
\HH(A, F)_\bullet = \mathrm{Bar}(A)_\bullet\otimes_{A\otimes A^\rev} F
= A^{\otimes\bullet}\otimes F
\]
regarded as a simplicial object $\Delta^{\mathrm{op}}\to\mathcal R$.
If $\mathcal R$ admits geometric realizations and the tensor product
preserves geometric realizations in each variable, the Hochschild
homology of $F$ exists and can be calculated via the Hochschild complex:
\[
\Tr(A, F) = |\HH(A, F)_\bullet|.
\]

An important special case of categorical trace occurs when the algebra
$A$ is equipped with an endomorphism $\phi$. Let $F$ be an
$A$-$A$-bimodule, then the $\phi$-twisted bimodule $F$, which we
denote by $^\phi F$, is the bimodule with the same right $A$-action
but with the left $A$-action precomposed with $\phi$. We write
$\Tr(A, \phi)$ for the Hochschild homology of $F = {}^\phi A$.

When $\mathcal R=\Lincatlam$, the $(\infty, 2)$ category of
$\Lambda$-linear $\infty$-categories, we refer to $\Tr(A, F)$ (or
$\Tr(A, \phi)$) as the categorical trace.

\makepoint[Abstract sheaf theory.]
Let $\mbf{C}$ be an $\infty$-category that admits finite limits and
finite coproducts, and let $\mathrm{pt}$ denote its final object.
The category $\mbf{C}$ will play the role of the category of geometric
objects.

Let $\cV, \cH$ be two weakly stable class of morphisms in $\mbf{C}$
(i.e., they contain all isomorphisms in $\mbf{C}$, are stable
under equivalences of morphisms, and are stable under base
change and composition). Let
$\Corr(\mbf{C})_{\mathrm{V}; \mathrm{H}}$ denote the category
of correspondences, as defined in \cite[\S~8..1.1]{Zhu25}.
We write $\mathrm{All}$ for the class of all morphisms in $\mbf{C}$,
and simply write $\Corr(\mbf{C})$ for
$\Corr(\mbf{C})_{\mathrm{All}; \mathrm{All}}$.

\begin{rmk}
\label{vr-hr}
There is a subclass of morphisms $\mathrm{VR}\subset\mathrm{V}$
(resp. $\mathrm{HR}\subset\mathrm{H}$) which are called
\emph{vertically} (resp. \emph{horizontally}) \emph{right
adjointable}. For their definitions see \cite[\S~8.2.2]{Zhu25}. We do
not recall the definitions here, but we will indicate the classes
$\mathrm{VR}$ and $\mathrm{HR}$ explicitly in our applications.
\end{rmk}

An abstract sheaf theory with coefficient in $\Lambda$ of
$\mbf{C}$ is a lax symmetric monoidal functor
\[
\shD: \mathrm{Corr}(\mathbf{C})_{\mathrm{V}; \mathrm{H}}\to
\Lincatlam.
\]
The functoriality of a sheaf theory encodes a rich structure. To begin,
we denote by $g^\star:\shD(X)\to\shD(Y)$ the functor corresponding to
a horizontal morphism $X\xleftarrow{g}Y\xrightarrow{\id}Y$, and we
denote by $f_\dagger:\shD(X)\to\shD(Y)$ the functor corresponding to
a vertical morphism $X\xleftarrow{\id}X\xrightarrow{f}Y$. A
general correspondence $X\xleftarrow{g}Z\xrightarrow{f}Y$ is
sent to a functor isomorphic to $f_\dagger\circ g^\star$. Moreover, the
sheaf theory encodes familiar structural results such as the base change
theorem and the projection formula. For further discussion, see
\cite[\S~8.2.1]{Zhu25}.

\makepoint[Convolution pattern.]
Let $X, Y\in\mathbf{C}$ and $f:X\to Y$. Then $X\times_Y X$ carries a
natural algebra structure in $\mathrm{Corr}(\mathbf{C})$ with the
multiplication and the unit maps given by
\[
\begin{tikzcd}[column sep=huge]
X \times_Y X \times_Y X
  \arrow[d, "\mathrm{id} \times f \times \mathrm{id}"]
  \arrow[r, "\mathrm{id} \times \Delta_X \times \mathrm{id}"]
  & (X \times_Y X) \times (X \times_Y X) \\
X \times_Y X
\end{tikzcd},
\qquad
\begin{tikzcd}
X \arrow[r] \arrow[d, "\Delta_{X/Y}"]
& \mathrm{pt} \\
X \times_Y X
\end{tikzcd}.
\]
Let $g_1: Z\to Y$ be another morphism, then $X\times_Y Z$ is a left
$(X\times_Y X)$-module with the action map given by
\[
\begin{tikzcd}[column sep=huge]
X \times_Y X \times_Y Z
  \arrow[d, "\mathrm{id}_X \times f \times \mathrm{id}_Z"]
  \arrow[r, "\mathrm{id}_X \times \Delta_X \times \mathrm{id}_Z"]
  & (X \times_Y X) \times (X \times_Y Z) \\
X \times_Y Z
\end{tikzcd}.
\]
Similarly, let $g_2: Z\to Y$ be another morphism, then
$X\times_{Y, g_1} Z\times_{g_2, Y} X$ is an $(X\times_Y X)$-bimodule.

Suppose there is an automorphisms $\sigma_X:X\to X$ and an
endomorphism $\sigma_Y:Y\to Y$ such that $f\circ\sigma_X =
\sigma_Y\circ f$. We will abuse notation and denote both maps by
$\sigma$ if it is clear from context. Let us take $Z = Y, g_1 = f$ and
$g_2 = f\circ \sigma$ in the discussion above.
Then the $(X\times_Y X)$-bimodule $X\times_Y Z\times_Y X$ is canonically
isomorphic to the $\sigma$-twisted module ${}^{\sigma}(X\times_Y X)$
by sending $(x, z, x')\in X\times_Y Z\times_Y X$ to $(\sigma(x), x')$.

Let us fix a sheaf theory $\shD:\mathrm{Corr}(\mathbf{C})_{\mathrm{V};
\mathrm{H}}\to\Lincatlam$ in the sequel.
Applying $\shD$ to our convolution pattern yields
an algebra obejct $\shD(X\times_Y X)$, a $\shD(X\times_Y X)$-bimodule
$\shD(X\times_Y Z\times_Y X)$, and an isomorphism
$\shD(X\times_Y Z\times_Y X) \cong {}^\sigma\shD(X\times_Y X)$. 

Consider the diagram
\begin{equation}
\label{pull-push}
\begin{tikzcd}[column sep=huge]
X \times_{Y\times Y} Z 
  \arrow[r, "\delta_0 = (\Delta_X \times \mathrm{id}_Z)"]
  \arrow[d, "q = (f \times \mathrm{id}_Z)"']
& (X \times X) \times_{Y\times Y} Z \\
Y \times_{Y\times Y} Z. &
\end{tikzcd}
\end{equation}
We have $(X \times X) \times_{Y\times Y} Z = X\times_Y Z\times_Y X$,
and $Y\times_{Y\times Y} Z = \Lsig Y$. The diagram hence
induces a functor $q_\dagger\circ\delta_0^\star:
\shD(X\times_Y Z\times_Y X)\to\shD(\Lsig Y)$, where $\Lsig Y$ is the
fixed point stack. The following proposition is
\cite[Corollary~8.60]{Zhu25}.

\begin{prop}
\label{cor860}
Assume $f:X\to Y\in\mathrm{VR}$ and
$\Delta_X:X\to X\times X\in\mathrm{HR}$, then there is a canonical
factorization
\[
\begin{tikzcd}[row sep=large, column sep=large]
\mathcal{D}(X \times_Y X)
  \arrow[r, "\delta_0^\star"]
  \arrow[d]
& \mathcal{D}(X \times_{X \times X} (X \times_Y X))
  \arrow[d, "q_{\dagger}"] \\
\operatorname{Tr}_{\mathrm{geo}}(\mathcal{D}(X \times_Y X), \sigma)
  \arrow[r]
& \mathcal{D}(\mathcal{L}_{\sigma} Y),
\end{tikzcd}
\]
where the lower horizontal arrow is fully faithful and its essential image
is generated under colimits by the image of $q_\dagger\circ\delta_0^\star$.
See Section~\ref{tr-geom} for the explanation of the notation
$\Tr_\mathrm{geo}$.
\end{prop}

\begin{definition}
Let $\mathrm{Eproet}$ denote the class of essentially
pro-\'etale morphisms in $\mathrm{AlgSp}^{\mathrm{perf}}_{\kappa_F}$,
i.e.\ those $f : X \to Y$ that can be written as $f : X \to X' \to Y$
with $X \to X'$ pro-\'etale and $X' \to Y$ perfectly finitely
presented (pfp).

We denote by $\mathrm{IndEproet}$ the class of morphisms
$f : X \to Y$ such that for every map $S \to Y$ with $S \in
\Sch^{\perf}_{\kappa_F}$, the pullback $f_S : X_S \to S$ admits a
presentation as a filtered colimit
\[
    X_S = \varinjlim_{i \in I} X_i
\]
with transition maps $X_i \to X_j$ being pfp closed immersions of
algebraic spaces and with each $f_i : X_i \to S$ belonging to
$\mathrm{Eproet}$.
\end{definition}

We will apply Proposition \ref{cor860} to the following situations.

\begin{eg}
\label{eg-shv}
$\mathbf{C} = \PreStk^\perf_{\kappa_F}$, with sheaf
theory
\[
\Shv(-, \Lambda): \Corr(\PreStk^\perf_{\kappa_F})
_{\IndEproet; \All} \to\Lincatlam,
\]
$\VR$ the class of ind-pfp proper morphisms, and $\HR$ the class of
representable pseudo cohomological pro-smooth morphisms.
A correspondence $X\xleftarrow{g}Z\xrightarrow{f}Y$ is associated to
a functor isomorphic to $g_!\circ f^*$.
$X = \bb{B}\LnT, Y = \bb{B}LT, f:X\to Y$ is induced by the inclusion
$\LnT\to LT$, $\sigma_X$ and $\sigma_Y$ are the Frobenius
automorphisms of $X$ and $Y$ respectively. (See
\cite[Proposition~10.97, Corollary~10.102]{Zhu25} for the definition
of this sheaf theory and the classes $\VR$ and $\HR$).
\end{eg}

\begin{eg}
\label{eg-coh}
$\mathbf{C} = \IndArStk^\mathrm{aft}_\Zl$ is the category
of ind-Artin stacks almost of finite presentation over $\Zl$, with
sheaf theory
\[
\IndCoh: \Corr(\IndArStk^\mathrm{aft}_\Zl)
_{\All; \mathrm{ftor}} \to\Lincatlam,
\]
$\VR$ the class of representable proper morphisms, and $\HR$ the
class of morphisms representable of finite tor-dimension
($\,\mathrm{ftor}$).
A correspondence $X\xleftarrow{g}Z\xrightarrow{f}Y$ is associated to
a functor isomorphic to $g_*\circ f^!$.
$X = Y = \X, f:X\to Y$ is the identity, $\sigma_X$
is the Frobenius automorphism of $X$. (See
\cite[Theorem~9.32]{Zhu25}).
\end{eg}

\refmakepoint
\label{tr-geom}
Let $X_1 = X\times_Y X$ and $Q = X\times_Y Z\times_Y X$.
We do not recall the definition of the geometric trace
$\Tr_\mathrm{geo}$ from \cite[\S~8.3.1]{Zhu25}. Instead, it is clear
that both Example \ref{eg-shv} and Example \ref{eg-coh} satisfy
the following conditions:
\begin{enumerate}[label=(\roman*)]
    \item $X\to Y$ is in VR.
    \item The diagonal $\Delta_X:X\to X\times X$ is in HR.
    \item The exterior tensor product
    \[
    \boxtimes_{\shD(\mathrm{pt})}:
    \shD(X_1)\otimes_{\shD(\mathrm{pt})}\shD(X_1)\to
    \shD(X_1\times X_1)
    \]
    is an equivalence.
    \item The exterior tensor product
    \[
    \boxtimes_{\shD(\mathrm{pt})}:
    \shD(X_1)\otimes_{\shD(\mathrm{pt})}\shD(Q)\to
    \shD(X_1\times Q)
    \]
    is an equivalence.
\end{enumerate}
Hence by \cite[Proposition 8.71]{Zhu25}, the comparison map
\[
\Tr(\shD(X_1), \shD(Q))\to
\Tr_\mathrm{geo}(\shD(X_1), \shD(Q))
\]
is an equivalence. This allows us to use Proposition \ref{cor860} with
the ordinary categorical trace in lieu of the geometric trace.

We first apply Proposition \ref{cor860} to Example \ref{eg-shv} to get
the following result.

\begin{cor}
\label{tr-shv}
Let $X = \bb BL^\n T$, $Y = \bb BLT$, and let $f:X\to Y$ be the
morphism induced by the inclusion $\LNT\to LT$. The diagram
(\ref{pull-push}) becomes
\begin{equation}
\begin{tikzcd}
    \dfrac{LT}{\Ad_\sigma L^\n T} \ar[r, "\delta_0"] \ar[d, "q"] &
    \doublequotientLT \\
    \dfrac{LT}{\Ad_\sigma LT},
\end{tikzcd}
\end{equation}
where $\delta_0$ is the identity on objects and sends $t \in
L^\n T$ to $(t, -\sigma t)$ on morphisms, while $q$ is the identity
on objects and the inclusion $L^\n T\to LT$ on morphisms. There
is a canonical factorization
\begin{equation}
\label{tr-series}
\begin{tikzcd}
    \Shv\!\Big(\doublequotientLT\Big) \ar[d]
    \ar[rd, start anchor={[yshift=-1em]east}, "q_!\circ\delta_0^*"]\\
    \Tr\Big(\Shv\!\Big(\doublequotientLT\Big), \sigma\Big)
    \ar[r, "\Psi_n"'] & \Shv(\isoc_T),
\end{tikzcd}
\end{equation}
where $\Psi_n$ is fully faithful.
\end{cor}

\refmakepoint \label{point-def-psi}
Let $s:*\to\frac{LT}{\Adsig LT}$ be a point, and form the pullback
diagram
\[
\begin{tikzcd}
    \dfrac{LT}{L^\n T} \ar[r, "s'"] \ar[d, "q'"]&
    \dfrac{LT}{\Ad_\sigma L^\n T} \ar[r, "\delta_0"] \ar[d, "q"]&
    \doublequotientLT \\
    * \ar[r, "s"]&
    \dfrac{LT}{\Ad_\sigma LT},
\end{tikzcd}
\]
where $s'$ is given by $t\mapsto t - \sigma t$ on objects and
identity on morphisms. We have
\begin{equation*}
    s^*\circ q_*\circ\delta_0^*\,\mathcal L
    = q'_*\circ s'^*\circ\delta_0^*\,\mathcal L.
\end{equation*}
Since $\frac{L^+T}{L^\n T}\cong T\times\bb{A}^d$ as
an algebraic scheme for some integer $d$, two observations follow.
First, the functors $q_!$ and $q_*$ differ only by a cohomological
shift. Second, the image of $q_*\circ\delta_0^*$ is independent
of $n$. Therefore, taking colimit of the diagram (\ref{tr-series})
yields the commutative diagram
\begin{equation}
\label{def-psi}
\begin{tikzcd}
    \Shv(LT) \ar[d]
    \ar[rd, start anchor={[yshift=-0.5em]east}, "Q"]\\
    \Tr(\Shv(LT), \sigma)
    \ar[r, "\Psi"'] & \Shv(\isoc_T),
\end{tikzcd}
\end{equation}
where $\Psi$ is fully faithful.

Now we apply Proposition \ref{cor860} to Example \ref{eg-coh}.
\begin{cor}
    \label{cor-def-phi}
    Let $X=\X$, then $\Lsig X = \Locc$. Let $i:\Locc\to\X$
    be the canonical morphism of the fixed-point stack.
    There is a canonical factorization
    \begin{equation}
    \label{def-phi}
    \begin{tikzcd}
        \IndCoh(\X)\ar[d]\ar[rd, start anchor={[yshift=1ex]}, "i^*"]\\
        \Tr(\IndCoh(\X), \sigma) \ar[r, "\Phi"'] & \IndCoh(\Locc),
    \end{tikzcd}
    \end{equation}
    where $\Phi$ is fully faithful, and the essential image of $i^*$
    generates $\IndCoh_\Nilp(\Locc)$.
\end{cor}
\begin{proof}
    The morphism $i$ is the composition of a closed embedding with a
    smooth morphism, therefore the functors $i^!$ and $i^*$ differs
    only by a cohomological shift.
    
    From Corollary~\ref{component-in-hat-h} we know that
    $\Sing(\X)$ splits as $\X\times\B(\mathfrak t^*)_I$.
    The image of $\Sing(i): \Sing(\X)\times_\X \Locc\to\Sing(\Locc)$
    is $\Locc\times\B(\mathfrak t^*)_\Gamma$. By
    \cite[Proposition~9.66]{Zhu25}, the essential image of $i^*$
    generates $\IndCoh_{\Im(\Sing(i))}(\Locc) = \IndCoh_\Nilp(\Locc).$
\end{proof}

\section{Geometric Local Langlands Correspondence for Tori}
\label{geom-llc}
\noindent
In this section, let $\Lambda = \Flb$ or $\Qlb$. By abuse of notation,
we use $\Locc$ and $\X$ to denote their base changes
$\Locc\otimes\Lambda$ and $\X\otimes\Lambda$.

The second main result of this paper is the following.

\begin{thm}
    \label{second-main}
    (i) Let $T$ be an algebraic torus over $F$. There exists a
    fully-faithful, $t$-exact, monoidal functor
    \[
        \Ch: \IndCoh(\X) \to \Shv(LT)
    .\]
    The essential image of $\Ch$ is compactly generated by the
    translations of all character sheaves on $L^+T$ to $LT$.
    
    (ii) Both $\IndCoh(\X)$ and $\Shv(LT)$ carry a canonical Frobenius
    structure. Let $\bb T$ denote the functor induced by $\Ch$:
    \[
    \bb T: \Tr(\IndCoh(\X), \sigma)\to \Tr(\Shv(LT), \sigma),
    \]
    and let $\bb L: \IndCoh_\Nilp(\Locc)\to\Shv(\isoc_T)$ be the
    categorical local Langlands correspondence under the
    identification $\QCoh(\Tor)\cong\Shv(\isoc_T)$. Then the following
    diagram commutes:
    \begin{equation}
    \begin{tikzcd}
        \label{decategorification}
        \Tr(\IndCoh(\X), \sigma)
        \ar[r, "\bb T"] \ar[d, "\Phi"]
        & \Tr(\Shv(LT), \sigma) \ar[d, "\Psi"] \\
        \IndCoh_\Nilp(\Locc) \ar[r, "\bb{L}"] & \Shv(\isoc_T).
    \end{tikzcd}
    \end{equation}
    (For definition of $\Phi$ and $\Psi$, see
    (\ref{def-psi}) and (\ref{def-phi}).) Furthermore, $\Phi$ and
    $\Psi$ are equivalences of categories.
\end{thm}

\makepoint[Construction of $\Ch$.]
We define the functor
$\Ch:\IndCoh(\X)\to\Shv(LT)$ as follows.

Using the split short exact sequence (\ref{ses3}), we have a
decomposition
\[
	\IndCoh(\X) = \bigoplus_{\alpha\in \hat L_I}
	\IndCoh^\alpha(\X)
.\] 
On the other hand, since
\[
	1\to L^+ T \to LT \to \hat L_I\to 1
,\] 
we have another decomposition
\[
	\Shv(LT) = \bigoplus_{\beta\in \hat L_I}
	\Shv(L^+T)
.\] 
The index set of the two decomposition are canonically identified.
However, for our purpose, the functor sends the $\alpha$-component
of ind-coherent sheaves to the $(-\alpha)$-component of
$\Lambda$-sheaves. Now we only have to define $\Ch$ on the
0-component.

Given that $\Coh^0(\X) = \Coh(\rfsc{\pi_1(L^+T), \Gm})$,
it is enough to define the functor
\[
    \Ch^0:\Coh(\rfsc{\pi_1(L^+T), \Gm})\to\Shv(L^+T),
\]
and naturally extend it to a functor between ind-completions
of both sides. Recall the decomposition from (\ref{decomp_X}):
\[
    \rfsc{\piLTvalue, \Gm}
    = \bigsqcup_{x\in\Xi} \mathcal{R}_{\hat J, \Gm}.
\]
For a component indexed by $x\in\Xi$, according to Proposition
\ref{embedding}, coherent sheaves on
$\mathcal{R}_{\hat J, \Gm}$ are given by
\[
    \Coh(\mathcal{R}_{\hat J, \Gm})
    = \underset{Z\subset \hat H}\colim \Coh(Z),
\]
where the colimit runs through all $Z$ finite over $\Zl$ such that
$Z(\Flb)\subset\hat H(\Flb)^p$. Each coherent sheaf $\mathcal F_x$ can
be regarded as a $\hat J$-module $V_{\mathcal F}$ that is finite over
$\Zl$.

On the other hand, each closed point
$x\in\Xi = \rfsc{(\hat L\otimes P)^I, \Gm}$
determines a rank-one $(\hat L\otimes P)^I$-module, which we
denote by $L_x$. We define the functor $\Ch^0$ as
\begin{align*}
        \Ch^0: \Coh(\rfsc{\piLTvalue, \Gm}) &\to
        \Rep\!\left(
        (\hat L\otimes I^\tame_E)^I \times
        (\hat L\otimes P)^I\right)\\
        \mathcal F_x &\mapsto V_{\mathcal F}\boxtimes L_x.
\end{align*}

\begin{proof}[Proof of Theorem \ref{second-main} (i)]
    By construction, $\Ch^0$ is $t$-exact, and hence so is $\Ch$.

    To show $\Ch$ is fully faithful, it is enough to show $\Ch^0$ is.
    Recall that $\Coh^0(\X)$ is generated by $\mathcal O_x$ for all
    closed points $x\in\X$. It suffices to check $\Ch^0$ is fully
    faithful on these generators.
    Let $x_1, x_2\in\X$ be closed points.
    By Corollary~\ref{component-in-hat-h}, the component in which
    $x_i$ lies is the formal \nhd in $\hat H = \mathcal R_{J, \Gm}$
    for $i = 1, 2$. Using the Koszul resolution of $\Lambda[X_1,
    \cdots, X_r]$, $r = \rank \hat L$, we get
    \[
        \Hom(\mathcal O_{x_1}, \mathcal O_{x_2}) =
        \begin{cases}
            0, & x_1\neq x_2,\\
            \bigwedge(L\otimes\Lambda[-1]),
            & x_1=x_2.
        \end{cases}
    \]
    Let $\mathcal L_i = \Ch^0(\mathcal O_{x_i})$ for $i = 1, 2$. We also
    have
    \[
        \Hom(\mathcal L_1, \mathcal L_2) =
        \begin{cases}
            0, & x_1\neq x_2,\\
            \bigwedge(L\otimes\Lambda[-1](-1))
            & x_1=x_2.
        \end{cases}
    \]

    Since both $\X$ and $LT$ carry abelian group structures,
    $\IndCoh(\X)$ and $\Shv(LT)$ carry natural monoidal
    structures by $*$-pushforward along the multiplication map.
    To show $\Ch$ is monoidal, we first decompose both categories into
    blocks
    \begin{align*}
        \IndCoh(\X) &= \bigoplus_{\alpha\in \hat L_I}
        \IndCoh^\alpha\left(\X\right),\\
        \Shv(LT) &= \bigoplus_{\beta\in \hat L_I}
        \Shv(L^+ T)
    \end{align*}
    where $\Ch$ respects both direct sums (with a twist $\beta = -\alpha$).
    If $\mathcal F \in \IndCoh^{\alpha_1}(\X)$ and
    $\mathcal G \in \IndCoh^{\alpha_2}(\X)$, we have
    $\mathcal F \star \mathcal G \in \IndCoh^{\alpha_1 + \alpha_2}(\X)$.
    Similar fact holds for $\Shv(LT)$, so it remains to show
    $\Ch^0$ is monoidal on the $0$-component.

    Since $\Coh^0(\X)$ is generated by all skyscraper sheaves
    $\Lambda_\chi$ indexed by characters $\chi:\cHom(\piLT, \lx)^\sigma$,
    it suffices to check $\Ch$ is monoidal on these generators.
    It is clear from our construction that for such $\chi$, $\Ch(\chi)$
    is nothing but the character sheaf associated to the character
    $\chi$. We hence have
    \[
    \Ch(\Lambda_{\chi_1} \star \Lambda_{\chi_2})
    = \Ch(\Lambda_{\chi_1 + \chi_2})
    = \Ch(\Lambda_{\chi_1}) \star \Ch(\Lambda_{\chi_2}).
    \]
\end{proof}

\makepoint[Decategorification.]
We move on to prove part (ii) of Theorem \ref{second-main}.
Consider the following diagram:
\begin{equation}
\label{two-squares}
\begin{tikzcd}[column sep=large]
    \IndCoh(\X) \ar[r, "\Ch"] \ar[d, "\mathrm{tr}"]
    \ar[dd, bend right=80, shift right=0, "i^*", swap]
    & \Shv(LT) \ar[d]
    \ar[dd, bend left=70, shift left=0, "Q"] \\
    \Tr(\IndCoh(\X), \sigma)
    \ar[r, "\bb T"] \ar[d, "\Phi"]
    & \Tr(\Shv(LT), \sigma) \ar[d, "\Psi"] \\
    \IndCoh_\Nilp(\Locc) \ar[r, "\bb{L}"] & \Shv(\isoc_T).
\end{tikzcd}
\end{equation}

The upper square is commutative by functorality of the categorical
trace construction.

By Corollary~\ref{cor-def-phi}, $\Phi$ is an equivalence of categories.
Since $\Psi$ is fully-faithful, it suffices to show that the lower
square commutes in order to deduce that $\Psi$ is an equivalence.

Since the essential image of $i^*$ generates
$\IndCoh_\Nilp(\Locc)$, to show the lower square commutes, it is
enough to show that the lower square commutes for the essential image
of $\mathrm{tr}$. Furthermore, by Corollary \ref{generator-for-coh},
$\mathcal O_x$ generates $\IndCoh(\X)$ for all closed points $x\in\X$.
If $x\neq\sigma x$, then $i^*\mathcal O_x = 0$. Thus it remains to
check
\[
        \mathbb{L}\circ i^*\mathcal O_x\cong
        Q\circ\Ch(\mathcal O_x)
\]
for every $x = \sigma x$.

Under Corollary~\ref{image-of-i} and Lemma~\ref{invariant-char},
let a $\sigma$-invariant closed point $x\in\X$ corresponds to a character
$\chi:T(F)_0\to\lx$. We let $I_\chi$ denote be the augmentation ideal
$\{t - \chi(t) | t\in T(F)_0\}$ in the group ring $\Lambda[T(F)]$ and
define $R_\chi = \Lambda[T(F)]/I_\chi$.

\begin{prop}
    \label{calculate-H}
    We have
    \[
            i^*\, \mathcal O_x
            = R_\chi \otimes \bigwedge (X_*(T)\otimes \Lambda[1]).
    \]
\end{prop}
\begin{proof}
\sloppy
The components of $\Locc$ are indexed by the set of representations
$\cHom(T(F)_0, \lx)$. The component $U\subset\Locc$ indexed by
$\chi$ maps to the closed point $x\in\cHom(\piLT, \lx)^\sigma$ and
$U = \Spec(R_\chi)$.

By Corollary~\ref{component-in-hat-h}, the component in which $x$ lies
can be embedded in some $\hat H = \mathcal R_{J, \Gm}$, and the
component itself is precisely the formal \nhd of $x$ in that torus.
Therefore,
\[
i^*\, \mathcal O_x = R_\chi \otimes^L_{\Lambda[J]} \Lambda
= R_\chi \otimes^L_{\Lambda[\hat L]} \Lambda
= R_\chi \otimes \bigwedge (X^*(T)\otimes \Lambda[-1]). \qedhere
\]
\end{proof}

Let $\mathcal L = \Ch(\mathcal O_x)$. By assumption, $\sigma^*\mathcal
L = \mathcal L$. We fix an $n$ such that $\chi$ factors through
$\frac{L^+T}{\LNT}(\kappa_F)$. The next proposition
computes $Q(\mathcal L)$. 

\begin{prop}
    We have
    \[
    Q(\mathcal L) = 
    R_\chi \otimes \bigwedge (X^*(T)\otimes\Lambda[-1]).
    \]
\end{prop}
\begin{proof}
Let $s:*\to\frac{LT}{\Adsig LT}$ be a point.
We calculate $Q(\mathcal L) = q_*\circ\delta_0^*\,\mathcal L$ using
the pullback diagram:
\[
\begin{tikzcd}
    \dfrac{LT}{L^\n T} \ar[r, "s'"] \ar[d, "q'"]&
    \dfrac{LT}{\Ad_\sigma L^\n T} \ar[r, "\delta_0"] \ar[d, "q"]&
    \doublequotientLT \\
    * \ar[r, "s"]&
    \dfrac{LT}{\Ad_\sigma LT},
\end{tikzcd}
\]
where $s'$ is given by $t\mapsto t - \sigma t$ on objects and
identity on morphisms.

We first determine the components of $\tfrac{LT}{\LNT}$ on
which the pullback $s'^*\circ\delta_0^*\,\mathcal L$ is non-zero. The
components of $\frac{LT}{\Adsig \LNT}$ are indexed by the set
$X_*(T)_I$, and so is $\frac{LT}{\LNT}$. Clearly,
$\delta_0^*\,\mathcal L$ has the same support as $\mathcal L$, which
is on the central component. Let $[s] \in X_*(T)_I$ be the component
to which $s$ belongs, then $s'^*\circ\delta_0^*\,\mathcal L$ has support
on all those components $\alpha\in X_*(T)_I$ such that 
\begin{equation*}
[s] + \alpha - \sigma\alpha = 0.
\end{equation*}
Assume $\alpha_0$ satisfies the equation above, then
$s'^*\circ\delta_0^*\,\mathcal L$ has support on those components
indexed by $\alpha_0 + (X_*(T)_I)^\sigma$.

Next, we determine on which components of $\frac{LT}{\Adsig LT}$ the
sheaf $q_*\circ\delta_0^*\,\mathcal L$ is non-zero. The components of
$\frac{LT}{\Adsig LT}$ are indexed by $B(T) = X_*(T)_\Gamma$. The fact
that $[s] = \sigma\alpha - \alpha$ in the previous analysis implies
that $s$ lies in the central component. Therefore, to study
$q_*\circ\delta_0^*\,\mathcal L$, it is enough to take $s = 1$.

By our assumption, $\sigma^*\, \mathcal L\cong \mathcal L$, hence the
pullback $s'^*\circ\delta_0^*\,\mathcal L$ is the trivial local system
on components indexed by $(X_*(T)_I)^\sigma$.

Since $\frac{L^+T}{L^\n T}\cong T\times\bb{A}^d$ as
an algebraic scheme for some integer $d$, we have
\begin{align*}
    s^*\circ q_*\circ\delta_0^*\,\mathcal L
    &= q'_*\circ s'^*\circ\delta_0^*\,\mathcal L\\
    &= \bigoplus_{(X_*(T)_I)^\sigma}
       \bigwedge \left(X^*(T)\otimes\Lambda[-1](-1)\right).
\end{align*}
The sheaf $q_*\circ\delta_0^*\,\mathcal L$ is the above sheaf equipped
with a $T(F)$-action. The action of $T(F)_0$ is given by the character
$\chi$. The action of $T(F)/T(F)_0\cong (X_*(T)_I)^\sigma$ is by
translating the direct summands. Combining the two actions, we obtain
\begin{equation*}
    q_*\circ\delta_0^*\,\mathcal L
    = R_\chi \otimes \bigwedge (X^*(T)\otimes\Lambda[-1](-1)).
    \qedhere
\end{equation*}
\end{proof}

\appendix

\section{Lemma in Homological Algebra}
\label{LHS-lem}

\begin{lem}
Let $\Gamma$ be a finite group, and let $A$ be an abelian group
equipped with a $\Gamma$-action. Let $[\alpha]\in H^2(\Gamma, A)$
represent the group extension
\[
1\to A\to G\to \Gamma\to 1.
\]
Let $M$ be a $\Gamma$-module, or equivalently, a $G$-module on which
$A$ acts trivially. Then we have a commutative diagram
\[
\hspace{-8ex}
\begin{tikzcd}
& H_2(\Gamma, M) \ar[r, "d_2"] \ar[d, "\cup\alpha"]
& H_1(A, M)_\Gamma \ar[r] \ar[d, equal]
& H_1(G, M) \ar[r] \ar[d, "\mathrm{res}"]
& H_1(\Gamma, M) \ar[r] \ar[d, "\cup\alpha"] & 0\\
0 \ar[r] & \hat{H}^{-1}(\Gamma, M\otimes A) \ar[r] &
(M\otimes A)_\Gamma \ar[r] &
(M\otimes A)^\Gamma \ar[r] &
\hat{H}^0(\Gamma, M\otimes A) \ar[r] & 0.
\end{tikzcd}
\]
The first row of the diagram is the long exact sequence associated to the
Lyndon-Hochschild-Serre spectral sequence. The second row is the definition
of the Tate cohomology.
The first and last vertical map is the cup product with $\alpha$. The map 
$\mathrm{res}$ is the restriction map $H_1(G, M)\to H_1(A, M)^\Gamma$.

In particular, if the triple $(\Gamma, A, \alpha)$ satisfy the
Tate-Nakayama criterion and if $M$ is torsion free as an abelian group,
then
\[
H_1(G, M) \cong (M \otimes A)^\Gamma.
\]
\end{lem}
\begin{proof}
    The crux of the proof is the commutativity of the left-most
    square. The commutativity of the remaining ones are fairly
    standard and will be omitted.

    \makepoint\textbf{Notations.}
    We follow the notation of Atiyah and Wall \cite{cassels1987algebraic}.
    The homogeneous bar resolution of the $G$-module $\Z$ is given by
    \[
    C_n(G, \Z) = \bigoplus_{G^{n+1}}\Z(g_0, g_1, \cdots, g_n),
    \]
    which is the free abelian group generated by the basis $G^{n+1}$
    with the action of
    \[
    g\cdot (g_0, \cdots, g_n) = (g g_0, \cdots, g g_n).
    \]
    The differential is the usual
    \[
    d(g_0, \cdots, g_n)
    = \sum_{i = 0}^n (-1)^i (g_0, \cdots, \widehat{g_i}, \cdots, g_n).
    \]
    This complex is identified with the more familiar inhomogeneous bar resolution, by sending
    \[
    (g_0, g_0 g_1, \cdots, g_0 g_1 \cdots g_n)\mapsto (g_0|g_1|\cdots|g_n).
    \]
    
    Let $P_\bullet$ denote a $\Gamma$-resolution of $\mathbb{Z}$ by
    finitely-generated free $\Gamma$-modules, and let $P^* =
    \mathrm{Hom}(P_\bullet, \mathbb{Z})$ be its dual, so that we have exact
    sequences
    \begin{align*}
        \cdots \to P_1 \to P_0 \to \mathbb{Z} \to &0 \\
        &0 \to \mathbb{Z} \to P_0^* \to P_1^* \to \cdots.
    \end{align*}
    Let $P_{-n} = P_{n-1}^*$ and splicing the two sequences together, we
    obtain a doubly-infinite exact sequence
    \[
    L_\bullet: \quad \cdots \to P_1 \to P_0 \to P_{-1} \to P_{-2} \to \cdots
    \]
    The Tate groups are defined as the cohomology groups of
    $\Hom_\Gamma(L_\bullet, M)$ for any $\Gamma$-module $M$, i.e.
    \[
    H^i_T(\Gamma, M) = H^i(\Hom_\Gamma(L_\bullet, M)).
    \]
    The homology of $M$ agrees with the Tate cohomology groups
    of degree $\leq -2$ via isomorphisms
    \[
    P_\bullet\otimes_\Gamma M = (P_\bullet\otimes M)_\Gamma \xrightarrow{N}
    (P_\bullet\otimes M)^\Gamma \to (\Hom(P_\bullet^*, M))^\Gamma
    = \Hom_\Gamma(P_\bullet^*, M).
    \]

    \makepoint\textbf{Cup Product.}
    Let $\{w_\sigma|\sigma\in\Gamma\}$ be a system of representatives
    of $\Gamma$ in $G$. The class of the extension $\alpha$ is
    represented by the cocycle
    $\delta: \Gamma\times\Gamma\to A$ defined by
    \[
    w_\sigma w_\tau = \delta(\sigma, \tau)\, w_{\sigma\tau}.
    \]
    The 2-cocycle $\delta$ determines a homogeneous 2-cocycle
    $\delta' \in \Hom_\Gamma(C_2(\Gamma, \mathbb Z), A)$ characterized by
    \[
    \delta'(1, c_1, c_1 c_2) = \delta(c_1, c_2).
    \]

    Let $[\phi]\in H_2(\Gamma, M)$ be represented by a 2-cycle
    $\phi: \Gamma\times\Gamma\to M$. Then it is also represented by a 
    homogeneous 2-cycle $\phi'\in\Hom_\Gamma(C_2^*(\Gamma, \Z), M)$
    such that
    \[
    \phi'(1, c_1^*, (c_1 c_2)^*) = \phi(c_1, c_2).
    \]
    The cup product $[\phi]\cup[\delta]\in H_0(\Gamma, M)$ is defined
    to be the class of
    $\phi'\cup\delta'\in\Hom_\Gamma(C_0^*(\Gamma, \Z), M\otimes A)$,
    which is the homogeneous cycle defined by
    \begin{align*}
        \phi'\cup\delta'(c^*) &= \sum_{c_1, c_2} \phi'(c^*, c_1^*, c_2^*)
        \otimes\delta'(c_2, c_1, c)\\
        &= \sum_{c_1, c_2} \phi'(c^*, c_1^*, (c_1 c_2)^*)
        \otimes\delta'(c_1 c_2, c_1, c).
    \end{align*}
    Its image in $(M\otimes A)_\Gamma$ is
    \begin{align*}
        \phi'\cup\delta'(1) &= \sum_{c_1, c_2} \phi'(1, c_1^*, (c_1 c_2)^*)
        \otimes\delta'(c_1 c_2, c_1, 1)\\
        &= \sum_{c_1, c_2} \phi(c_1, c_2)\otimes
        c_1 c_2\cdot\delta'(1, c_2^{-1}, c_2^{-1} c_1^{-1})\\
        &= \sum_{c_1, c_2} \phi(c_1, c_2)\otimes
        c_1 c_2\cdot\delta(c_2^{-1}, c_1^{-1}).
    \end{align*}

    \makepoint\textbf{Spectral sequence.} For a $G$-module $M$, we now
    construct the double complex giving rise to the
    Lydon--Hochschild--Serre spectral sequence. Let $P_\bullet =
    C_\bullet(\Gamma, \Z)$ and $Q_\bullet = C_\bullet(G, \Z)$. The
    differential of the two chain complexes are denoted by $d_1$ and
    $d_0$ respectively. We define the double complex
    \[
    E_{ij} = P_i\otimes Q_j.
    \]
    The total complex $\mathrm{Tot}(E_{ij})$ forms a $G$-resolution of
    $\Z$, hence the spectral sequence associated with $E_{ij}\otimes_G
    M$ calculates the group homology of $M$. Note that $P_\bullet$ is
    also a $G$-module. We have
    \[
    E_{ij}^M := E_{ij}\otimes_{\,G} M = (P_i\otimes Q_j\otimes M)_G = \left(
    P_i\otimes_{\,\Gamma} (Q_j\otimes_A M)\right).
    \]
    Since $Q_j\otimes_A (-)$ calculates $A$-homology and $P_i\otimes_{\,\Gamma}
    (-)$ calculates $\Gamma$-homology, the resulting spectral sequence of this
    double complex is the Lydon-Hochschild-Serre spectral sequence.

    Next, We relate the $A$-homology calculated via $Q_\bullet$ to
    that calculated via $C_\bullet(A, \Z)$. We choose a projection
    $f: G\to A$ such that
    \begin{align*}
        f(w_\gamma) &= 1, && \text{for all } \gamma\in\Gamma,\\
        f(a g) &= a f(g), && \text{for all } a\in A, g\in G.
    \end{align*}
    This choice induce a morphism between $A$-complexes $C_\bullet(G,
    \Z)\to C_\bullet(A, \Z)$. In particular, $C_1(G, \Z)\otimes_A M\to
    C_1(A, M)$ is given by
    \begin{equation}
        \label{send-to-am}
        (g_1, g_2)\otimes m \mapsto (f(g_2) - f(g_1))\otimes m.
    \end{equation}

    \makepoint\textbf{Differentials.} Recall that in the sepctral
    sequence, the differential $d_2: H_2(\Gamma, M) \to H_1(A,
    M)_\Gamma$ is defined as follows. Let $\phi\in E_{2, 0}^M$ be a
    cycle and assume there exists $\psi\in E_{1, 1}^M$ such that
    \begin{equation*}
    d_0\psi + d_1\phi = 0,
    \end{equation*}
    then
    \[
    d_2([\phi]) = d_1([\psi]).
    \]

    From now on, we use $c$, $c_i$ or Greek letters $\sigma, \tau$ to
    denote elements of $\Gamma$, $g, g_i$ or $h$ to denote
    elements of $G$. We note that in $E_{i, j}^M$, the following identity
    holds
    \begin{align*}
        &\quad\,(c_0, c_1, \cdots, c_i)\otimes (g_0, \cdots, g_j)\otimes m\\
        &= (1, c_0^{-1} c_1, \cdots, c_0^{-1} c_i)
        \otimes (w_{c_0}^{-1} g_0, \cdots, w_{c_0}^{-1} g_j)
        \otimes c_0^{-1} m.
    \end{align*}
    Consequently, every cycle $\phi\in E_{i, j}^M$ can be expressed as
    \[
        \phi = \sum_{c_1, \cdots, c_i, g_1, \cdots, g_j}
        (1, c_1, c_1 c_2, \cdots, c_1 c_2 \cdots c_i)
        \otimes (g_0, \cdots, g_j)\otimes
        \phi(c_1, c_2, \cdots, c_i, g_1, \cdots, g_j).
    \]
    In general, this expression is not unique. However, when $j = 0$,
    it becomes unique up to the following identity
    \[
    (1, c_1, c_1 c_2, \cdots, c_1 c_2 \cdots c_i)\otimes(g)\otimes m
    = (1, c_1, c_1 c_2, \cdots, c_1 c_2 \cdots c_i)
    \otimes(a g)\otimes m
    \]
    for any $a\in A$. Equivalently, if we define
    \[
    \tphi(c_1, \cdots, c_i, \sigma) = 
    \sum_{\bar g = \sigma} \phi(c_1, \cdots, c_i, g),
    \]
    then $\phi$ admits a unique representation
    \[
    \phi = \sum_{c_1, \cdots, c_i, \sigma}
    (1, c_1, c_1 c_2, \cdots, c_1 c_2 \cdots c_i)\otimes
    (w_\sigma)\otimes \tphi(c_1, \cdots, c_i, \sigma).
    \]
    
    Now we proceed to compute $d_2$. The next three formulae are
    straightforward.
    \begin{align*}
        \phantom{=}& d_1\Big((1, c_1, c_1 c_2)\otimes(g)\otimes m\Big)\\
        =& \left( c_1 (1, c_2) - (1, c_1 c_2) + (1, c_1)\right)\otimes(g)\otimes m\\
        =& (1, c_2)\otimes (w_{c_1}^{-1} g) \otimes c_1^{-1} m
        -(1, c_1, c_2)\otimes(g)\otimes m + (1, c_1)\otimes(g)\otimes m.\\[1em]
        \phantom{=}& d_0\,(1, c)\otimes(g_0, g_1)\otimes m\\
        =& (1, c)\otimes(g_0)\otimes m - (1, c)\otimes(g_1)\otimes m.\\[1em]
        \phantom{=}& d_1\,(1, c)\otimes (g_1, g_2)\otimes m\\
        =& (w_c^{-1} g_1, w_c^{-1} g_2)\otimes c^{-1} m - (g_1, g_2)\otimes m\\
        \mapsto& c^{-1} m \otimes \left(f(w_c^{-1} g_2) - f(w_c^{-1} g_1)\right)
        - m \otimes \left(f(g_2) - f(g_1)\right)\\
        =& c^{-1} m \otimes \left(\delta(c^{-1}, g_2) - \delta(c^{-1}, g_1)\right)
        - m \otimes \left(f(g_2) - f(g_1)\right).
    \end{align*}
    In the second to last line, we send the element to $M\otimes A$
    using Equation (\ref{send-to-am}). Applying the formulas above, we
    compute the differentials explicitly as follows.
    \begin{align*}
        d_1\,\tphi(\gamma, \sigma)
        &= \sum_{c_2 = \gamma} c_1^{-1}\tphi(c_1, c_2, c_1 \sigma)
        - \sum_{c_1 c_2 = \gamma} \tphi(c_1, c_2, \sigma)
        + \sum_{c_1 = \gamma} \tphi(c_1, c_2, \sigma)\\
        &= \sum_{c} c^{-1} \tphi(c, \gamma, c\sigma)
        - \sum_{c} \tphi(c, c^{-1} \gamma, \sigma)
        + \sum_{c} \tphi(\gamma, c, \sigma),\\[1em]
        d_0\,\tpsi(\gamma, \sigma)
        &= \sum_{\overline{h} = \sigma,\,g} \psi(\gamma, h, g) - \sum_g \psi(\gamma, g, h).\\[1em]
        d_1\,\psi &= \sum_{c, g_1, g_2} c^{-1} \psi(c, g_1, g_2)\otimes
        (\delta(c^{-1}, g_2)-\delta(c^{-1}, g_1))\\
        &\phantom{=} - \sum_{c, g_1, g_2} \psi(c, g_1, g_2)\otimes (f(g_2)-f(g_1)).
    \end{align*}
    
    Finally, assume $d_1\,\phi + d_0\,\psi = 0$. Then we obtain
    \begin{align*}
        d_2([\phi]) &= d_1\,\psi\\
        &= \sum_{\gamma, g_1, g_2}
        \gamma^{-1}\psi(\gamma, g_1, g_2)\otimes (\delta(\gamma^{-1}, g_2)-\delta(\gamma^{-1}, g_1))\\
        &= \sum_{\gamma, g_1, g_2} 
        \big(\gamma\cdot\psi(\gamma^{-1}, g_2, g_1) - \gamma\cdot\psi(\gamma^{-1}, g_1, g_2)\big)\otimes \delta(\gamma, g_1)\\
        &= \,\sum_{\gamma, \sigma} - \gamma\cdot d_0\,\tpsi(\gamma^{-1}, \sigma)\otimes \delta(\gamma, \sigma)\\
        &= \,\sum_{\gamma, \sigma} \gamma\cdot d_1\,\tphi(\gamma^{-1}, \sigma)\otimes \delta(\gamma, \sigma)\\
        &= \sum_{\gamma, \sigma, c} 
        \big(c^{-1}\tphi(c, \gamma^{-1}, c \sigma) - \tphi(c, c^{-1}\gamma^{-1}, \sigma)
        + \tphi(\gamma^{-1}, c, \sigma)\big) \otimes \delta(\gamma, \sigma)
        \\
        &= \sum_{c_1, c_2, \sigma} \tphi(c_1, c_2, \sigma)\otimes \big(
        c_1 c_2 \cdot \delta(c_2^{-1}, c_1^{-1} \sigma) - c_1 c_2 \cdot \delta(c_2^{-1} c_1^{-1}, \sigma)
        + c_1 \delta(c_1^{-1}, \sigma) \big)\\
        &= \sum_{c_1, c_2, \sigma} \tphi(c_1, c_2, \sigma)\otimes c_1 c_2\cdot \delta(c_2^{-1}, c_1^{-1})\\
        &= [\phi]\cup\alpha. \qedhere
    \end{align*}
\end{proof}

%% file: LLT.bib
@article{Zhu21,
      title={Coherent sheaves on the stack of Langlands parameters}, 
      author={Xinwen Zhu},
      year={2021},
      eprint={2008.02998},
      archivePrefix={arXiv},
      primaryClass={math.AG},
      url={https://arxiv.org/abs/2008.02998}, 
}

@article{BB07,
  author  = {Braverman, Alexander and Bezrukavnikov, Roman},
  title   = {Geometric {Langlands} correspondence for {D}-modules in prime characteristic: the {GL(n)} case},
  journal = {Pure and Applied Mathematics Quarterly},
  year    = {2007},
  volume  = {3},
  number  = {1},
  pages   = {153--179},
  note    = {Special Issue: In honor of Robert D. MacPherson. Part 3},
  doi     = {10.4310/PAMQ.2007.v3.n1.a5}
}

@article{CZ14,
  author  = {Chen, Tsao-Hsien and Zhu, Xinwen},
  title   = {Geometric {Langlands} in prime characteristic},
  journal = {Compositio Mathematica},
  year    = {2017},
  volume  = {153},
  number  = {2},
  pages   = {395--452},
  doi     = {10.1112/S0010437X16008113}
}

@article{Kot85,
  author  = {Kottwitz, Robert E.},
  title   = {Isocrystals with additional structure},
  journal = {Compositio Mathematica},
  year    = {1985},
  volume  = {56},
  number  = {2},
  pages   = {201--220}
}

@misc{Kot14,
      title={B(G) for all local and global fields}, 
      author={Robert Kottwitz},
      year={2014},
      eprint={1401.5728},
      archivePrefix={arXiv},
      primaryClass={math.RT},
      url={https://arxiv.org/abs/1401.5728}, 
}

@article{Zhu14,
  author  = {Zhu, Xinwen},
  title   = {Affine Grassmannians and the geometric {Satake} in mixed characteristic},
  journal = {Annals of Mathematics},
  year    = {2017},
  volume  = {185},
  number  = {2},
  pages   = {403--492}
}

@article{Lan97,
  title={Representations of abelian algebraic groups},
  author={Langlands, Robert P.},
  journal={Pacific Journal of Mathematics, Special Issue},
  year={1997},
  volume={181},
  pages={231--250}
}

@article{Zou24,
      title={The categorical form of Fargues' conjecture for tori}, 
      author={Konrad Zou},
      year={2024},
      eprint={2202.13238},
      archivePrefix={arXiv},
      primaryClass={math.RT}
}

@book{SGA4Expose18,
    title={{SGA 4 Expos{\'e} XVIII}},
    author={Deligne, P.}
}

@article{DW23,
      title={Character sheaves on tori over local fields}, 
      author={Tanmay Deshpande and Saniya Wagh},
      year={2023},
      eprint={2304.06622},
      archivePrefix={arXiv},
      primaryClass={math.RT},
      url={https://arxiv.org/abs/2304.06622}, 
}

@article{Ser61,
     author = {Serre, Jean-Pierre},
     title = {Sur les corps locaux \`a corps r\'esiduel alg\'ebriquement clos},
     journal = {Bulletin de la Soci\'et\'e Math\'ematique de France},
     pages = {105--154},
     publisher = {Soci\'et\'e math\'ematique de France},
     volume = {89},
     year = {1961},
     doi = {10.24033/bsmf.1562},
     mrnumber = {26 #103},
     zbl = {0166.31103},
     language = {fr}
}

@article{Zhu25,
      title={Tame categorical local Langlands correspondence}, 
      author={Xinwen Zhu},
      year={2025},
      eprint={2504.07482},
      archivePrefix={arXiv},
      primaryClass={math.RT},
      url={https://arxiv.org/abs/2504.07482}, 
}

@book{cassels1987algebraic,
  editor    = {J.~W.~S. Cassels and A.~Fr{\"o}hlich},
  title     = {{Algebraic Number Theory: Proceedings of an Instructional Conference}},
  publisher = {Academic Press},
  address   = {London},
  year      = {1987},
  edition   = {Reprint of the 1967 volume},
  pages     = {366},
  isbn      = {0--12--163251--2},
}

@article{Breen,
  author  = {Lawrence Breen},
  title   = {Un th\'eor\`eme d'annulation pour certains $\mathrm{Ext}^i$ de faisceaux ab\'eliens},
  journal = {Annales scientifiques de l'\'Ecole normale sup\'erieure, 4\textsuperscript{e} s\'erie},
  volume  = {8},
  number  = {3},
  pages   = {339--352},
  year    = {1975},
  doi     = {10.24033/asens.1291},
}
